\documentclass[10pt]{article}
\usepackage{epsfig}
\title{An iterative thresholding algorithm for linear inverse problems
with a sparsity constraint}
\author{Ingrid Daubechies \and Michel Defrise \and
Christine De Mol}
\usepackage{amssymb,amsmath,amsfonts}

\newtheorem{theorem}{Theorem}[section]
\newtheorem{proposition}[theorem]{Proposition}
\newtheorem{lemma}[theorem]{Lemma}
\newtheorem{remark}[theorem]{Remark}
\newtheorem{corollary}[theorem]{Corollary}

\newcommand{\QED}{\hfill \raisebox{-2pt}{\rule{5.6pt}{8pt}\rule{4pt}{0pt}}%
  \smallskip\par}

\def\R{\mathbb{R}}
\def\C{\mathbb{C}}
\def\N{\mathbb{N}}
\def\Z{\mathbb{Z}}
\def\S{\mathbf{S}}
\def\T{\mathbf{T}}
\def\I{\sl {I}}
\def\A{\mathbf{A}}
\def\cB{\mathcal{B}}
\def\cH{\mathcal{H}}

\def\mR{\mbox{R}}
\def\mN{\mbox{N}}
\def\K{K^*}
\def\wf{\widetilde{f}}
\def\fs{f^{\dagger}}
\def\OBJ{\mbox{\tiny\texttt{OBJECT}}}
\def\IM{\mbox{\tiny\texttt{IMAGE}}}
\def\SUR{\mbox{\tiny{\it{SUR}}}}
\def\fg{f_{\gamma}}
\def\ag{w_{\gamma}}
\def\vpg{\varphi_{\gamma}}
\def\g{\gamma}
\newcommand\eref[1]{(\ref{#1})}
\def\Vvert{\vert\!\vert\!\vert}
\def\VVert{[\!|\!|\!]}
\begin{document}
\maketitle

\begin{abstract}

We consider linear inverse problems where the solution is assumed to have
a sparse expansion on an arbitrary pre--assigned orthonormal basis.
We prove that replacing the usual
quadratic regularizing penalties by weighted $\ell^p$-
penalties on the coefficients of such expansions, with $1 \leq p \leq 2$,
still regularizes the problem. If $p < 2$, regularized solutions of such
$l^p$-penalized problems will have sparser expansions, with respect to the
basis under consideration. To compute the corresponding regularized solutions
we propose an iterative algorithm which amounts to a Landweber iteration
with thresholding (or nonlinear shrinkage) applied at each iteration step.
We prove that this algorithm converges in norm.  We
also review some potential applications of this method.

\end{abstract}

\section{Introduction}
\subsection{Linear inverse problems}
In many practical problems in the sciences and applied sciences, the features
of most interest cannot be observed directly, but have to be inferred
from other, observable quantities. In the simplest
approximation, which works surprisingly well in a wide range of cases
and often suffices,
there is a {\em linear} relationship between the features of interest
and the observed quantities. If we model the {\em object} (the traditional
shorthand
for the features of interest) by a function $f$, and the derived quantities or
{\em image} by another function $h$, we can cast the problem of inferring
$f$ from $h$
as a {\em linear inverse problem}, the task of which is to
solve the equation
\begin{equation*}
Kf=h~.
\end{equation*}
This equation and the task of solving it make sense only when placed in an
appropriate
framework.
In this paper, we shall assume that $f$ and $h$ belong to appropriate function
spaces, typically Banach or Hilbert spaces,
$f \in \cB_{\OBJ}$,
$h \in \cB_{\IM}$, and that $K$ is a bounded  operator from
the space  $\cB_{\OBJ}$ to
$ \cB_{\IM}$. The choice of the spaces must be
appropriate for describing real-life situations.

The observations or {\em data}, which we shall model by yet another
function, $g$, are typically not
exactly equal to the image $h=Kf$, but rather to a distortion of $h$. This
distortion
is often modeled by an {\em additive noise} or {\em error term} $e$, i.e.
\begin{equation*}
g=h+e=Kf+e ~.
\end{equation*}
Moreover, one typically assumes that the ``size'' of the noise can be measured
by its $L^2$-norm, $\|e\|=\left( \int_{\Omega} |e|^2\right)^{1/2}$ if $e$
is a function on $\Omega$. (In a finite-dimensional situation, one uses
$\|e\|=
\left(\sum_{n=1}^N |e_n|^2 \right)^{1/2}$ instead.)
Our only ``handle'' on the image $h$ is thus via the observed $g$, and
we typically have little information on $e=g-h$ beyond an upper bound on
its
$L^2$-norm $\|e\|$.
(We have here implicitly placed ourselves
in the ``deterministic setting'' customary to the
Inverse Problems community. In the
stochastic setting more familiar to statisticians, one assumes instead
a bound on the variance
of the components of $e$.)
Therefore it is customary to take
$\cB_{\IM}=L^2(\Omega)$; even if the ``true images'' $h$ (i.e. the images
$Kf$ of the
possible objects $f$) lie in a much smaller space, we can only know them up
to some (hopefully) small $L^2$-distance.

We shall consider in this paper
a whole family of possible choices for $\cB_{\OBJ}$, but we shall
always assume that these spaces are subspaces of a basic Hilbert space
$\cH$ (often
an $L^2$-space as well), and that $K$ is a bounded operator from $\cH$ to
$L^2(\Omega)$.
In many applications, $K$ is an integral operator with a kernel
representing the response of the imaging device; in the special case where
this linear device is
translation-invariant, $K$ reduces to a convolution operator.

To find an estimate of $f$ from the observed $g$, one can minimize
the {\it discrepancy} $\Delta(f)$,
\begin{equation*}
\Delta(f)=\| Kf-g \|^2\ ;
\end{equation*}
functions that minimize $\Delta(f)$ are called {\em pseudo-solutions}
of the inverse problem. If the operator $K$ has a trivial null-space,
i.e. if $\mN(K)=\{f \in \cH ; Kf=0\}=\{0\}$, there is a unique minimizer, given
by
$\widetilde{f}=(K^*K)^{-1}K^*g$, where $K^*$ is the adjoint operator. If
$\mN(K)\neq
\{0\}$ it is customary to choose, among the set of pseudo-solutions, the unique
element
$f^\dagger$ of minimal norm,
i.e. $f^\dagger = \mbox{arg-min}\{\|f\|; f \mbox{ minimizes }
\Delta(f)\}$.
This function belongs to $\mN(K)^\perp$ and is called the {\it generalized
solution}
 of
the inverse problem. In this case the map $K^\dagger:
g \mapsto f^\dagger$
is called the {\it
generalized inverse } of $K$. Even when $K^*K$ is not invertible,
$K^\dagger g $ is
well-defined for all $g$ such that $K^*g \in \mR(K^*K)$.
 However, the generalized inverse operator may be
unbounded (for so-called {\it ill-posed problems}) or have a very large norm
(for {\it ill-conditioned problems}). In such instances, it has to be replaced
by bounded approximants or approximants with smaller norm, so
that numerically stable solutions can be defined and used as meaningful
approximations of the true solution corresponding to the exact data.
This is the issue of {\it regularization}.

\subsection{Regularization by imposing additional quadratic constraints}

The definition of a pseudo-solution (or even, if one considers equivalence
classes modulo $\mN(K)$, of a generalized solution) makes use of
the inverse of the operator $K^*K$; this inverse is well defined on
the range $\mR(K^*)$ of $\K $ when $K^*K$ is a strictly
positive operator, i.e. when its spectrum is bounded below away from zero.
When the spectrum of $K^*K$ is not bounded below by a strictly positive
constant, $\mR(K^*K)$ is not closed, and not all elements of $\mR(\K )$ lie
in $\mR(\K K)$. In this case there is no guarantee
that $K^*g
\in \mR(\K K)$; even if $K^*g$ belongs to $\mR(\K K)$, the unboundedness of
$(\K K)^{-1}$ can cause severe numerical instabilities unless
additional measures are taken.

This blowup or these numerical instabilities are ``unphysical'', in the sense
that we know a priori that the true object would not have had a huge norm in
$\cH$, or other characteristics exhibited by the unconstrained ``solutions''.
A standard procedure to avoid these instabilities or to {\em regularize}
the inverse
problem is to modify the functional to be minimized, so that it
incorporates not
only  the discrepancy, but also the a priori knowledge one may have about the
solution.  For instance, if it is known that the object is of limited ``size''
in $\cH$, i.e. if
$\|f\|_{_{\cH}} \leq \rho$ , then the functional to be minimized can be chosen
as
\begin{equation*}
\Delta(f) + \mu \|f\|_{_{\cH}}^2 = \|Kf-g\|_{L^2(\Omega)}^2 +\mu
\|f\|_{_{\cH}}^2
\end{equation*}
where $\mu$ is some positive constant called the {\it regularization
parameter}. The minimizer is given by
\begin{equation}
\label{tikh}
f_\mu=(\K K + \mu I)^{-1} K^* g ~.
\end{equation}
where $I$ denotes the identity operator.
The constant $\mu$ can then be chosen appropriately,
depending on the application.
If $K$ is a
compact operator, with singular value decomposition given by
$Kf =\sum_{k=1}^{\infty} \sigma_k \left<f,v_k\right> u_k ~$, where $(u_k)_{k
\in \N}$ and $(v_k)_{k \in \N}$ are the orthonormal bases of eigenvectors of
$KK^*$ and $\K K$, respectively, with corresponding eigenvalues $\sigma_k^2$,
then \eref{tikh} can be rewritten as
\begin{equation}
\label{tikh-b}
f_\mu= \sum_{k=1}^{\infty} \frac{\sigma_k}{\sigma_k^2 + \mu}
\left<g, u_k\right> v_k ~.
\end{equation}
This formula shows explicitly how this regularization method reduces the
importance
of the eigenmodes of $\K K$ with small eigenvalues, which otherwise (if
$\mu =0$)
lead to instabilities. If an estimate of the ``noise'' is known, i.e. if
we know a priori that $g=Kf+e$ with $\|e\| \leq \epsilon$, then one finds
from \eref{tikh-b} that
\begin{equation*}
\| f - f_\mu \| \leq \left\Vert \sum_{k=1}^{\infty}
\frac{ \mu \left< f,v_k\right> }{\sigma_k^2+\mu}\; v_k \right\Vert +
\left\Vert \sum_{k=1}^{\infty} \frac{\sigma_k}{\sigma_k^2+\mu}
\left< e,u_k\right> v_k \right\Vert
\leq \Gamma(\mu) + \frac{\epsilon}{\sqrt\mu} ~,
\end{equation*}
where $\Gamma(\mu) \rightarrow 0$ as $\mu \rightarrow 0$. This means that
$\mu$ can be chosen appropriately, in an $\epsilon$-dependent way,
so that the error in estimation
$\|f - f_\mu \|$ converges to zero when $\epsilon$ (the estimation of
the noise level) shrinks to zero. This feature of the method, usually called
{\em stability}, is one that is required for any regularization method.
It is similar to requiring that a statistical estimator is consistent.

Note that the ``regularized estimate'' $f_\mu$ of \eref{tikh-b} is
linear
in $g$. This means that we have effectively defined a linear regularized
estimation operator that
is especially well adapted to the properties of the operator $K$; however,
it proceeds with a ``one method fits all'' strategy, independent of the data.
This may not
always be the best approach. For instance, if $\cH$ is an $L^2$-space
itself, and $K$
is an integral operator, the functions $u_k$ and $v_k$ are typically fairly
smooth;
if on the other hand the objects $f$ are
likely to have local singularities or discontinuities,
an approximation of type \eref{tikh-b}
(effectively limiting the estimates $f_\mu$
to expansions in the first $N ~v_k$, with $N$ determined by, say, $\sigma_k^2<
\mu/100$ for $k>N$) will of necessity be a smoothened version of $f$, without
sharp features.

Other classical regularization methods with
quadratic constraints may use quadratic Sobolev
norms, involving a few derivatives, as the ``penalty'' term added to the
discrepancy. This introduces a penalization of the
highly oscillating components, which are often the
most sensitive to noise. This method is especially
easy to use in the case where $K$ is a
convolution operator, diagonal in the Fourier domain.
In this case the regularization
produces a smooth cut-off on the highest Fourier components,
independently of the data.
 This works well
for recovering smooth objects which have their relevant structure contained in
the lower part of the spectrum and which have spectral content homogeneously
distributed across the space or time domain.
However, the Fourier domain is clearly not the appropriate representation
for expressing smoothness properties
of objects that are either spatially inhomogeneous, with
varying ``local frequency'' content, and/or present some
discontinuities, because
the frequency
cut-off implies that the resolution with which the fine details of the
solution can be stably retrieved is necessarily limited; it also
implies that the achievable
resolution is essentially the same at all points (see e.g. the book
\cite{Ber98}
for an extensive discussion of these topics).

\subsection{Regularization by non-quadratic constraints that promote sparsity}

The problems with the standard regularization methods described above
are well known and several approaches have been proposed for dealing
with them.
We propose in this paper a regularization method that, like the classical
methods just discussed, minimizes a functional obtained by adding a
penalization term
to the discrepancy; typically this penalization term will {\em not} be
quadratic,
but rather a weighted $\ell^p$-norm of the coefficients of $f$ with respect to
a particular orthonormal basis in $\cH$, with $1\leq p\leq 2$. More precisely,
given an orthonormal basis $\left(\vpg\right)_{\gamma\in \Gamma}$
of $\cH$, and given a sequence of strictly positive weights
${\mathbf w}= (\ag )_{\gamma\in \Gamma}$,
we define the functional $\Phi_{\mathbf{w},p}$ by
\begin{equation}
\label{funct-gen}
\Phi_{\mathbf{w},p}(f)= \Delta(f)
+\sum_{\gamma \in \Gamma} \ag\; |\left<f,\vpg\right>|^p
= \|Kf-g\|^2 +\sum_{\gamma \in \Gamma} \ag |\left<f,\vpg\right>|^p ~.
\end{equation}

For the special case $p=2$ and $\ag=\mu$ for all $\gamma \in \Gamma$
(we shall write this as $\mathbf{w}= \mu \mathbf{w}_{_0}$, where
$\mathbf{w}_{_0}$ is the sequence with all entries equal to 1),
this reduces
to the functional \eref{tikh}.  If we consider the family of functionals
$\Phi_{\mu\mathbf{w}_{_0},p}(f)$, keeping the weights
fixed at $\mu$, but decreasing $p$ from 2 to 1, we gradually
{\em increase} the penalization
on ``small'' coefficients (those with $|\left<f,\vpg\right>| <1$)
while simultaneously {\em decreasing} the
penalization on ``large coefficients'' (for which
$|\left<f,\vpg\right>|>1$).  As far as the
penalization term is concerned,  we are thus putting a lesser penalty on
functions $f$
with large but few components with respect to the basis
$\left(\vpg\right)_{\gamma\in \Gamma}$, and a higher penalty on
sums of many small components, when compared to the classical method of
\eref{tikh}.
This effect is the more pronounced the
smaller $p$ is. By taking $p <2$, and especially
for the limit value $p=1$,
the proposed minimization procedure thus promotes
sparsity of the expansion of $f$ with respect
to the $\vpg$.
(We shall not consider $p < 1$ here,
because then the functional
ceases to be convex.)

The bulk of this paper deals with  algorithms to obtain minimizers $f^*$
for the
functional
\eref{funct-gen}, for general operators $K$. In the special case where $K$
happens
to be diagonal in the $\vpg$--basis,
$K \vpg= \kappa_{\g} \vpg$,
the analysis is easy and straightforward.
Introducing the shorthand notation $\fg$ for $\left<f,\vpg\right>$ and
$g_{\gamma}$ for
$\left<g,\vpg\right>$, we have then
$$
\Phi_{\mathbf{w},p}(f)=
\sum_{\gamma \in \Gamma} \left[ |\kappa_{\g} \fg-g_{\gamma}|^2
+ \ag |\fg|^p\right] ~.
$$
The minimization problem thus uncouples into a family of 1-dimensional
minimizations, and is
easily solved. Of particular interest is the almost trivial case where
(i) $K$ is the identity operator,
(ii) $\mathbf{w}=\mu \mathbf{w}_{_0}$ and (iii) $p=1$,
which corresponds to the practical situation where
the data $g$ are equal to a noisy version of $f$ itself, and we want to remove
the noise
(as much as possible), i.e. we wish to {\em denoise} $g$. In this case the
minimizing $f^\star$ is given by
\begin{equation}
\label{simple}
f^\star = \sum_{\gamma \in \Gamma} f^\star_{\g} \vpg
= \sum_{\gamma \in \Gamma} S_{\mu}(g_{\gamma}) \vpg~,
\end{equation}
where $S_{\mu}$ is the (nonlinear) thresholding operator from $\R$ to $\R$
defined by
\begin{equation}
\label{stau}
S_{\mu}(x)= \left\{
 \begin{array}{ccl} x +\mu/2 ~&~ \mbox{if} ~& x \leq - \mu/2  \\
0 ~&~ \mbox{if} ~& |x| < \mu/2
\\ x- \mu/2 ~&~ \mbox{if} ~& x \geq  \mu/2.
\end{array} \right.
\end{equation}
(We shall revisit the derivation of \eref{simple} below. For simplicity,
we are assuming that all functions are real-valued. If the $\fg$ are complex,
a derivation similar
to that of \eref{simple} then leads to a complex thresholding operator, which
is defined as $S_{\mu}(r e^{i\theta})= S_{\mu}(r) e^{i\theta}$;
see Remark \ref{2-5} below.)

In more general cases, especially when $K$ is not diagonal with respect
to the $\vpg$--basis, it is not as straightforward to minimize
\eref{funct-gen}.

An approach that promotes sparsity with respect to a particular basis
makes sense only if we know that the objects $f$ that we want
to reconstruct do indeed have a sparse expansion with respect to this basis.
In the next subsection we list some situations in which this is the
case and to which the algorithm that we propose in this paper
could be applied.

\subsection{Possible applications for sparsity-promoting constraints}
\subsection*{1.4.1 Sparse wavelet expansions}

This is the application that was the primary motivation for this paper.
Wavelets provide orthonormal bases of $L^2(\R^d)$ with localization
in space and in scale; this makes them more suitable than e.g.
 Fourier expansions for an efficient
representation of functions that have space-varying smoothness properties.
Appendix \ref{WavBes} gives a very succinct overview of wavelets and their
link
with
a particular family of smoothness spaces, the Besov spaces. Essentially,
the Besov space $B^s_{p,q}(\R^d)$ is a space of functions on
$\R^d$ that ``have $s$ derivatives in $L^p(R^d)$''; the
index $q$ provides some extra fine-tuning. The
precise definition involves the moduli of continuity of the function,
defined by finite differencing, instead of derivatives, and combines
the behavior of these moduli at different scales. 
The Besov space $B^s_{p,q}(\R^d)$ is well-defined as
a complete metric space even if the indices $p,~q \in (0,\infty)$ are
$<1$, although it is no longer a Banach space in this case. 
Functions that are mostly smooth, but that have a few local
``irregularities'', nevertheless can still belong to a Besov space with
high smoothness index. For instance, the 1-dimensional function
$F(x)= \mbox{sign} (x)\; e^{-x^2}$ can belong to $B^s_{p,q}(\R)$ for
arbitrarily large $s$, provided $0<p<\left(s+\frac{1}{2}\right)^{-1}$. (Note that
this same example does not belong to any of the Sobolev spaces
$W^s_p(\R)$ with $s>0$, mainly because these can be defined only for
$p\geq 1$.) Wavelets provide unconditional bases for the Besov spaces,
and one can express whether or not a function $f$ on $\R^d$ belongs
to a Besov space by a fairly simple and completely explicit requirement
on the absolute values of the wavelet coefficients of $f$. 
This expression becomes particularly simple when $p=q$; as 
reviewed in Appendix \ref{WavBes},
$f \in B^s_{p,q}(\R)$ if and only if (see Appendix \ref{WavBes})
\begin{equation*}
\VVert f \VVert_{ _{s,p}} = \left( \sum_{\lambda \in \Lambda} 2^{\sigma p
|\lambda|} |\left<f, \Psi_{\lambda} \right> |^p \right) ^{1/p} < \infty~,
\end{equation*}
where $\sigma$ depends on $s,p$ and is defined by $\sigma =s +d \left(
\frac{1}{2}-\frac{1}{p} \right)$,
and where $|\lambda |$ stands for the scale of the wavelet
$\Psi_{\lambda}$. (The $\frac{1}{2}$ in the formula for $\sigma$ is
due to the choice of normalization
of the $\Psi_{\lambda}$, $\|\Psi_{\lambda}\|_{_{L^2}} =1$.)
For $p=q\geq 1$, $\VVert ~ \VVert_{ _{s,p}}$ is an equivalent norm to the
standard Besov norm on $B^s_{p,q}(\R^d)$; we shall restrict ourselves to 
this case in this paper.

It follows that minimizing
the variational functional for an inverse
problem with a Besov space prior falls exactly within the category of
problems studied in this paper: for such an inverse problem,
with operator $K$ and with
the a priori knowledge that the object lies in some $B^s_{p,p}$, it
is natural to define the
variational functional to be minimized
by
$$
\Delta(f)+ \VVert f \VVert _{ _{s,p}}^p = \|Kf-g\|^2 + \sum_{\lambda \in
\Lambda}
2^{\sigma p |\lambda|} |\left< f, \Psi_{\lambda} \right> |^p ~~,
$$
which is exactly of the type
$\Phi_{\mathbf{w},p}(f)$, as defined in \eref{funct-gen}.
For the case where $K$ is the identity operator, it was
noted already in \cite{Cha98}
that the wavelet-based algorithm for denoising
of data with a Besov prior, derived earlier in \cite{Don94},
amounts exactly to the minimization of
$\Phi_{\mu\mathbf{w}_{ _0},1}(f)$, where
$K$ is the identity operator and the $\vpg$--basis is a wavelet basis; the
denoised approximant given in \cite{Don94} then coincides
exactly with (\ref{simple}, \ref{stau}).

It should be noted that if $d > 1$, and if we are interested in functions that
are mostly smooth, with possible jump discontinuities
(or other ``irregularities'')
on smooth manifolds of dimension 1 or higher (i.e. not point irregularities),
then the Besov spaces do not constitute the optimal smoothness
space hierarchy. For
$d=2$, for instance, functions $f$ that are $C^{\infty}$ on the
square $[0,1]^2$, except on a finite set of smooth curves, belong to
$B^1_{1,1}([0,1]^2)$, but not to $B^s_{1,1}([0,1]^2)$ for $s>1$.
In order to obtain
more efficient (sparser) expansions of this type of functions, other expansions
have to be used, using e.g. ridgelets or curvelets (\cite{Don00},
\cite{Can00}).
One can then again use the approach in this paper, with respect to these
more adapted bases.

\subsection*{1.4.2 Other orthogonal expansions}

The framework of this paper applies to enforcing sparsity of the expansion of
the solution on any orthonormal basis. We provide here three examples which are
particularly relevant for applications, but this is of course not limitative.

The first example is the case where it is known a priori
that the object to be recovered is sparse in the Fourier domain, i.e. $f$
has only
a few nonzero Fourier components. It makes then sense to choose a standard
 Fourier
basis for the $\vpg$, and to apply the algorithms explained later in this
paper.
(They would have to be adapted to deal with complex functions, but
this is easily done; see Remark \ref{2-5} below.)
In the case where $K$ is the identity operator,
this is a classical problem, sometimes referred to as ``tracking sinusoids
drowned
in noise'', for which many other algorithms have been
developed.

For other applications, objects are naturally sparse in the original
(space or time) domain. Then our framework can be used again if we expand such
objects in a basis formed by the characteristic functions of pixels or voxels.
Once the inverse problem is discretized in pixel space, it is regularized by
penalizing the $l^p$-norm of the object with $1 \le p \le 2$. Possible
applications include the restoration of astronomical images with scattered
stars
on a flat background. Objects formed by a few spikes are also typical of some
inverse
problems arising in spectroscopy or in particle sizing. In medical
imaging, $l^p$-norm
penalization with $p$ larger than but close to $1$ has been used e.g.
for the restoration of tiny blood vessels \cite{Li02}.

The third example refers to the case where $K$ is compact and the use of SVD
expansions is a viable  computational approach, e.g. for the solution of
relatively  small-scale problems or
for operators that can be diagonalized in an analytic way. As already stressed
above, the linear regularization methods as given e.g. by \eref{tikh-b}
have the drawback that the penalization or cut-off eliminates the components
corresponding to the
smallest singular values, independently of the type of data. In some
instances,
the most
significant coefficients of the object may not correspond to the largest
singular
values; it may then happen that the object possesses
significant coefficient beyond
the cut-off imposed by linear methods. In order to avoid the elimination of
such
coefficients, it is preferable to use instead a
nonlinear regularization analogous to (\ref{simple}, \ref{stau}), with basis
functions $\vpg$ replaced
by the singular vectors $v_k$.
The theorems in this paper show that the {\it thresholded SVD expansion}
$$
f^* = \sum_{k=1}^{+\infty} S_{\mu/\sigma_k^2}
\left(\frac{\left< g,u_k\right>}{\sigma_k}
\right) v_k
~=~ \sum_{k=1}^{+\infty}\frac{1}{\sigma_k^2} ~ S_{\mu}
\left( \sigma_k \left< g,u_k\right>
\right) v_k~,
$$
which is the minimizer of the functional \eref{funct-gen}  with
$\mathbf{w}=\mu \mathbf{w}_{_0}$ and $p=1$,
provides a regularized solution that is better adapted to these
situations.

\subsection*{1.4.3 Frame expansions}

In a Hilbert space $\cH$, a frame $\{\psi_n\}_{n \in \N}$ is a
set of vectors for which there exist
constants $A, B >0$ so that, for all $v \in \cH$,
$$
B^{-1} \sum_{n\in \N} |\left<v,\psi_n\right>|^2 \leq \|v\|^2 \leq
A^{-1} \sum_{n \in \N}
|\left< v, \psi_n\right>|^2  ~~.\
$$
Frames always span the whole space $\cH$, but the frame vectors $\psi_n$
are
typically not linearly independent. Frames were first proposed
by Duffin and Schaeffer in \cite{DuSh52}; they are now used in a wide
range of
applications.
For particular choices of the frame vectors, the two frame bounds $A$
and
$B$ are equal; one has then, for all $ v \in \cH$,
\begin{equation}
\label{frame-1}
v = A^{-1} \sum_{n \in \N} \left<v, \psi_n \right> \psi_n ~.
\end{equation}
In this case, the frame is called {\em tight}.
An easy example of a frame is given by taking the union of two (or more)
different
orthonormal bases in $\cH$; these unions always constitute tight frames,
with $A=B$ equal to the number of orthonormal bases used in the union.

Frames are typically ``overcomplete'', i.e. they still span all of
$\cH$ even if some frame vectors are removed. It follows that, given
a vector $v$ in $\cH$, one can find many different sequences of
coefficients
such that
\begin{equation}
\label{frame-v}
v = \sum_{n \in \N} z_n \psi_n ~~.
\end{equation}
Among these sequences, some have special properties for which they
are preferred. There is, for instance,
a standard procedure to find the unique sequence
with minimal $\ell^2$-norm; if the frame is tight, then this
sequence is given by $z_n=A^{-1} \left<v,\psi_n\right>$,
as in \eref{frame-1}.

The problem of finding sequences
$\mathbf{z}=(z_n)_{n \in \N}$ that satisfy \eref{frame-v} can be
considered as an inverse problem. Let us define the operator $K$
from $\ell^2 (\N)$ to $\cH$ that maps a sequence
$\mbox{{\bf z}}= (z_n)_{n \in \N}$ to the element $K \mbox{{\bf z}}$ of $\cH$
by
\begin{equation*}
K \mbox{{\bf z}} = \sum_{n \in \N} z_n \psi_n ~.
\end{equation*}
Then solving \eref{frame-v} amounts to solving $K\mathbf{z}=v$. Note
that
this operator $K$ is associated with, but not identical to
what is often called the ``frame operator''. One has, for
$v \in \cH$,
$$
K K^* v = \sum_{n \in \N} \left<v, \psi_n \right> \psi_n~ ;
$$
for $\mathbf{z} \in \ell^2$, the sequence $K^*K\mathbf{z}$ is given by
$$
(K^*K \mathbf{z})_k = \sum_{l \in \N} z_l \left< \psi_l, \psi_k \right> ~.
$$
In this framework, the sequence $\mathbf{z}$ of minimum
$\ell^2$-norm that satisfies \eref{frame-v} is given simply by
$\mathbf{z}^{\dagger}= K^{\dagger}v$. The standard procedure in frame
lore
for the construction of $\mathbf{z}^{\dagger}$ can be rewritten as
$\mathbf{z}^{\dagger}=K^*(KK^*)^{-1}v$, so that
$K^{\dagger}=K^*(KK^*)^{-1}$
in this case. This last formula holds because this inverse problem is
in fact well-posed: even though $\mN(K) \neq \{0\}$, there is a gap
in the spectrum of $K^*K$ between the eigenvalue 0 and the remainder
of the spectrum, which is contained in the interval $[A,B]$; the
operator
$KK^*$ has its spectrum completely within $[A,B]$. In practice, one
always
works with frames for which the ratio $B/A$ is reasonably close to 1, so
that the problem is not only well-posed but also well-conditioned.

It is often of interest, however,
to find sequences
that are sparser than $\mathbf{z}^{\dagger}$.
For instance, one may know a priori that $v$ is a ``noisy'' version
of a linear combination of $\psi_n$ with a coefficient sequence of small
$\ell^1$-norm. In this case, it makes sense to determine a sequence
$\mathbf{z}_{\mu}$ that minimizes
\begin{equation}
\label{frame-3}
\|K\mathbf{z}-v\|^2_{\cH} + \mu \|{\mathbf{z}}\|_{\ell^1} ~,
\end{equation}
a problem that
falls exactly in the category of problems described in subsection 1.3.
Note that although the inverse problem for $K$ from $\ell^2(\N)$ to $\cH$
is well-defined, this need not be the case with the restriction
$K \big|_{\ell^1}$ from $\ell^1(\N)$ to $\cH$. One can indeed find
tight frames for which
$\sup \{ \|\mathbf{z}\|_{_{\ell^1}}~;~ \mathbf{z} \in \ell^1 \mbox{ and }
\|K\mathbf{z}\| \leq 1 \} = \infty$,
so that for arbitrarily large $R$ and arbitrarily small $\epsilon$,
one can find $\tilde{v} \in \cH$, $\tilde{\mathbf{z}} \in \ell^1$
with $\|\tilde{v}-K\tilde{\mathbf{z}}\| = \epsilon$, yet
$\inf \{ \|\mathbf{z}\|_{_{\ell^1}} ~; ~ \|\tilde{v}-K{\mathbf{z}} \|
\leq \epsilon/2 \} \geq R \|\tilde{\mathbf{z}}\|_{\ell^1}$.
In a noisy situation, it therefore may not make sense to search for the
sequence
with minimal $\ell^1$--norm that is ``closest'' to $v$; to find
an estimate of the $\ell^1$--sequences of which a given $v$ is known
to be a small perturbation, a better strategy is to compute the minimizer
$z_{\mu}$ of  \eref{frame-3}.

Minimizing the functional \eref{frame-3} as an approach to obtain
sequences that provide sparse approximations $K\mathbf{z}$
to $v$ was
proposed and applied to various problems by Chen, Donoho and Saunders
\cite{CDS01};
in the statistical literature, least-squares regression with
$\ell^1$-penalty is known as the ``lasso'' \cite{Tib96}.
The algorithm in this paper provides thus an alternative to linear and
quadratic programming techniques for these problems,
which all amount to minimizing  \eref{frame-3}.


\subsection{A summary of our approach}

Given an operator $K$ from $\cH$ to itself (or, more generally, from
$\cH$ to $\cH'$), and an orthonormal basis
$(\vpg)_{\gamma \in \Gamma}$, our goal is to find minimizing
$f^\star$ for the functionals $\Phi_{\mathbf{w},p}$ defined in section
1.3. The corresponding variational equations are
\begin{equation}
\label{variational}
\forall \gamma \in \Gamma ~ : ~
\left< \K Kf, \vpg \right> - \left< \K g , \vpg \right>
+ \frac{ \ag p}{2}~ | \left< f, \vpg\right> |^{p-1}
\mbox{sign}(\left< f, \vpg\right>) = 0 ~~.
\end{equation}
When $p \neq 2$ and $K$ is not diagonal in the
$\vpg$-basis, this gives a coupled system of
nonlinear equations for the $\left<f, \vpg \right>$,
which it is not immediately clear how to solve.
To bypass this problem, we
shall use a sequence of ``surrogate'' functionals that are each
easy to minimize,
and for which we expect, by a heuristic argument, that the successive
minimizers have our desired $f^\star$ as a limit.
These surrogate
functionals  are introduced in section 2 below. In section 3 we then show that
their successive minimizers do indeed converge to $f^\star$; we first
establish weak convergence, but conclude the section by proving that the
convergence also holds in norm. In section 4 we show that our proposed
iterative method is {\em stable}, in the sense given in subsection 1.2: if we
apply the algorithm to data that are a small perturbation of a ``true image''
$K f_o$,
then the algorithm will produce $f^\star$ that converge
to $f_o$ as the norm of the perturbation tends to zero.

\subsection{Related work}

Exploiting the sparsity of the expansion on a given basis
of an unknown signal, in order to assist in the estimation
or approximation of the signal from noisy data, is not a new
idea. The key role played by sparsity to achieve superresolution
in diffraction-limited imaging was already emphasized by Donoho
\cite{Don92} more than a decade ago. Since the seminal paper by Donoho
and Johnstone \cite{Don94}, the use of thresholding techniques for
sparsifying the wavelet expansions of noisy signals in order to remove the
noise (the so-called ``denoising'' problem) has been abundantly discussed in
the literature, mainly in a statistical framework (see e.g.
the book \cite{Mal98}).
Of particular importance for the background of this paper is the article by
Chambolle et al. \cite{Cha98}, which provides a variational formulation
for denoising, through the use of penalties on a Besov-norm of
the signal; this is the perspective adopted in the present paper.

Several attempts have been made to generalize the denoising framework
to solve inverse problems. To overcome the coupling problem stated in the
preceding subsection, a first approach is to construct wavelet- or
``wavelet-inspired'' bases that are in some sense adapted to the operator to
be inverted. The so-called Wavelet-Vaguelette
Decomposition (WVD) proposed by Donoho \cite{Don95}, as well as the twin
Vaguelette-Wavelet Decomposition method \cite{Abr98}, and also the
deconvolution in mirror wavelet bases \cite{KMR03, Mal98} can all be
viewed as examples of this strategy. For the inversion of the Radon transform,
Lee and Lucier \cite{Lee01} formulated a
generalization of the WVD decomposition that uses a variational
approach to set thresholding levels. A drawback
of these methods is that they are limited to special types of operators $K$
(essentially convolution-type operators under some additional
technical assumptions).

Other papers have explored the application of Galerkin-type methods to inverse
problems, using an appropriate but fixed wavelet basis \cite{Dic96, Lou97,
Coh02}. The underlying intuition is again that if the operator lends itself
to a
fairly  sparse representation in wavelets, e.g. if it is an operator of the
type
considered in \cite{Bey91}, and if the object is mostly smooth with some
singularities, then the inversion of the truncated operator will not be too
onerous, and the approximate representation of the object will do a good
job of
capturing the singularities. In \cite{Coh02}, the method is made
adaptive, so that
the finer-scale wavelets are used where lower scales indicate the presence of
singularities.

The mathematical framework in this paper has the advantage of not
pre-supposing any
particular properties for the operator $K$ (other than boundedness) or the
basis
$(\vpg)_{\gamma\in \Gamma}$ (other than its orthonormality). We prove,
in complete
generality, that generalizing Tikhonov's regularization method from the
$\ell^2$-penalty case to a $\ell^1$-penalty (or, more generally, a weighted
$\ell^p$-penalty with $1\leq p\leq 2$) provides a proper regularization
method for
ill-posed problems in a Hilbert space $\cal H$, with estimates that are
independent
of the dimension of $\cal H$ (and are thus valid for infinite-dimensional
separable $\cal H$). To our knowledge, this is the first proof of this fact.
Moreover, we derive a Landweber-type iterative algorithm that involves a
denoising procedure at each iteration step and provides a sequence of
approximations
converging in norm to the variational minimizer, with estimates of the rate of
convergence in particular cases. This algorithm was derived previously
in \cite{DeM02}, using,
as in this paper, a construction based on ``surrogate
functionals''. During the final editing of the present paper, our attention
was
drawn to the independent work by Figueiredo and Nowak \cite{Fig03, Nov01},
who, working in the different (finite-dimensional
and stochastic) framework of Maximum Penalized
Likelihood Estimation
for inverting a
convolution operator acting on objects that are sparse in the wavelet domain,
derive essentially the same iterative algorithm as in \cite{DeM02} and
this paper.

\section{An iterative algorithm through
surrogate functionals}
It is the combined presence of $\K Kf$ (which couples all the equations)
and the nonlinearity of the equations that makes the system
\eref{variational} unpleasant. For this reason, we borrow a technique
of optimization transfer (see e.g. \cite{Lan00} and \cite{DeP95}) and
construct surrogate functionals that effectively remove the term $\K Kf$. We
first pick a constant $C$ so that $ \|\K K \| < C$, and then
we define the functional $\Xi(f;a)=
C\|f-a\|^2-\|Kf-Ka\|^2$ which depends on an auxiliary element $a$ of
$\cH$.
Because $C\mbox{\I} - \K K$ is a strictly positive operator, $\Xi(f;a)$ is
strictly convex in $f$ for any choice of $a$. If $\|K\|<1$,
we are allowed to set $C=1$; for simplicity, we will restrict ourselves
to this case, without loss of generality since
$K$ can always be renormalized.
We then add  $\Xi(f ; a)$ to $\Phi_{\mathbf{w},p}(f)$ to form the following
``surrogate functional''
\begin{eqnarray}
\label{sur}
{\Phi}^{^{\SUR}}_{\mathbf{w},p}(f ; a)&=& \Phi_{\mathbf{w},p}(f)
- \|Kf-Ka\|^2 +  \|f-a\|^2 \nonumber \\
&=& \|Kf-g\|^2 + \sum_{\gamma \in \Gamma} \ag |\left< f,
\vpg \right>|^p -  \|Kf-Ka\|^2 +  \|f-a\|^2 \nonumber \\
&=& \|f\|^2 - 2\left<f, a+\K g  - \K K a \right> + \sum_{\gamma}
\ag | \left< f, \vpg \right>|^p + \|g\|^2 + \|a\|^2 - \| K a
\|^2\nonumber \\ &=& \sum_{\gamma} \left[ \fg^2 -2 \fg
\left(a +\K g  - \K K a \right)_{\gamma} + \ag |\fg|^p
\right] + \|g\|^2 + \|a\|^2 - \| K a \|^2
\end{eqnarray}
where we have again used the shorthand $v_{\gamma}$ for $\left<v,\varphi_
{\gamma} \right>$, and implicitly assumed
that we are dealing with real functions only.
Since $\Xi(f;a)$ is strictly convex in $f$,
${\Phi}^{^{\SUR}}_{\mathbf{w},p}(f ; a)$ is also strictly convex
in $f$, and has a unique minimizer for any choice of $a$. The advantage of
minimizing \eref{sur} in place of \eref{variational} is that the variational
equations for the $\fg$ decouple. We can then try to approach the minimizer
of $\Phi_{\mathbf{w},p}(f)$ by an iterative process which goes as follows:
starting from an arbitrarily chosen $f^0$, we determine the minimizer
$f^1$ of \eref{sur} for $a = f^0$; each successive iterate $f^n$ is then
the minimizer for $f$ of
the surrogate functional \eref{sur} anchored at the previous iterate, i.e. for
$a= f^{n-1}$. The iterative algorithm thus goes as follows
\begin{equation}
\label{iter}
f^0 {\mbox\ {\rm arbitrary}}\ ; f^n= \mbox{{\rm arg--min}}
\left({\Phi}^{^{\SUR}}_{\mathbf{w},p}(f ; f^{n-1})\right)\ \ n=1, 2,\dots
\end{equation}
To gain some insight into this iteration,
let us first focus on two special cases.

In the case where $\mathbf{w}=\mathbf{0}$ (i.e. the functional
$\Phi_{\mathbf{w},p}$ reduces to the discrepancy only), one needs to
minimize
$$
{\Phi}^{^{\SUR}}_{\mathbf{0},p}(f ; f^{n-1})=\|f\|^2-2\left<f,
f^{n-1} + \K (g - K f^{n-1}) \right> +\|g\|^2 + \|f^{n-1}\|^2 - \|Kf^{n-1}\|^2
~;
$$
this leads to
\begin{equation*}
f^{n} = f^{n-1} + \K (g - Kf^{n-1})~~.
\end{equation*}
 This is
nothing else than the so-called Landweber iterative method, the convergence of
which  to the (generalized) solution of
$Kf=g$ is well-known (\cite{Lad51}; see also \cite{Ber98},
\cite{Eng96}).

In the case where $\mathbf{w}=\mu \mathbf{w_0}$ and $p=2$,
the $n$-th surrogate functional reduces to
$$
{\Phi}^{^{\SUR}}_{\mathbf{w},2}(f ; f^{n-1})=(1+\mu)\;
\|f\|^2-2\left<f,f^{n-1} + \K (g - K f^{n-1}) \right> +\|g\|^2 +
\|f^{n-1}\|^2 -
\|Kf^{n-1}\|^2 ~;
$$
the minimizer is now
\begin{equation}
\label{DamLan}
f^n = \frac{1}{1+\mu} \left[ f^{n-1} +\K (g -K f^{n-1}) \right] ~,
\end{equation}
i.e. we obtain a damped or regularized Landweber iteration
(see e.g. \cite{Ber98}). The convergence of the function $f^n$ defined by
(\ref{DamLan}) follows immediately from the estimate
$\|f^{n+1}-f^n\| = (1+\mu)^{-1} \|(\mbox{\I}-\K K)(f^n-f^{n-1})\|
\leq (1+\mu)^{-1} \|f^n -f^{n-1}\|$, showing that we have a contractive
mapping,
even if $\mN(K) \neq \{0\}$.

In these two special cases we thus find that the $f^n$ converge as
$n \rightarrow \infty$. This permits one to hope that the $f^n$ will
converge for general $\mathbf{w}, p$ as well; whenever this is the
case the difference
$\|f^{n}-f^{n-1}\|^2 - \|K(f^{n}-f^{n-1})\|^2$ between
${\Phi}^{^{\SUR}}_{\mathbf{w},p}(f^n ; f^{n-1})$ and
$\Phi_{\mathbf{w},p}(f^n)$ tends to zero as $n \rightarrow \infty$,
suggesting that the minimizer
$f^n$ for the first functional could well tend to a minimizer
$f^\star$ of the second.
In section 3 we shall see that all this is more than a pipe-dream; i.e. we
shall prove that the $f^n$ do indeed converge to a minimizer  of
$\Phi_{\mathbf{w},p}$.

In the remainder of this section, we derive an explicit formula
for the computation of the successive
$f^n$. We first discuss the minimization of the functional
\eref{sur} for a generic $a \in \cH$.
As already noticed, the variational equations for the $\fg$
decouple. For $p>1$, the summand in \eref{sur} is differentiable in $\fg$, and
the minimization reduces to solving the variational equation
\begin{equation*}
2 \fg + p \, \ag \, \mbox{sign}(\fg) |\fg|^{p-1} = 2( a_{\gamma} + [\K (g-K
a)]_{\gamma}) ~;
\end{equation*}
since for any  $w \geq 0$ and any $p>1$,
the real function $F_{w,p}(x)=x+ {\frac{w p}{2}} ~ \mbox{sign}(x)|x|^{p-1}$ is
a one-to-one map
from $\mathbb{R}$ to itself,
we thus find that the minimizer of \eref{sur} satisfies
\begin{equation}
\label{solcomp-pneq1}
\fg= S_{\ag,p}( a_{\gamma} + [\K (g-K a)]_{\gamma}) ~,
\end{equation}
where $S_{w,p}$ is defined by
\begin{equation}
\label{S-pneq1}
S_{w,p}= \left( F_{w,p} \right)^{-1} ~, ~ \mbox{for } p>1.
\end{equation}

When $p=1$, the summand of  \eref{sur} is differentiable in $\fg$
only if $\fg \neq 0$; except at the point of non-differentiability,
the variational equation now reduces to
\begin{equation*}
2 \fg + \ag \,  \mbox{sign}(\fg) = 2 ( a_{\gamma} + [\K (g-K a)]_{\gamma})
~.
\end{equation*}
For $\fg>0$, this leads to $\fg= a_{\gamma} + [\K (g-K
a)]_{\gamma} -
\ag/2$; for consistency we must impose
in this case that $a_{\gamma} + [\K (g-K
a)]_{\gamma}
> \ag/2$. For $\fg <0$, we obtain
$\fg=a_{\gamma} + [\K (g-K
a)]_{\gamma}+\ag/2$, valid only when
$a_{\gamma} + [\K (g-K
a)]_{\gamma} < -\ag/2$. When
$a_{\gamma} + [\K (g-K
a)]_{\gamma}$ does not satisfy either of the two
conditions. i.e. when $|a_{\gamma} + [\K (g-K
a)]_{\gamma}| \leq \ag/2$,
we put $\fg =0$. Summarizing,
\begin{equation}
\label{solcomp-peq1}
\fg = S_{\ag,1}(a_{\gamma}+[ \K( g - K a)]_{\gamma}) ~,
\end{equation}
where the function $S_{w,1}$ from $\mathbb{R}$ to itself is defined by
\begin{equation}
\label{S-peq1}
S_{w,1}(x)=\left\{ \begin{array}{ccl} {x-w/2} & {\mbox{if}} & {x \geq w/2} \\
{0} & {\mbox{if}} & {|x| < w/2} \\ {x+w/2 }& {\mbox{if}} & {x \leq -w/2  ~~.}
\end{array} \right.
\end{equation}
(Note that this is the same nonlinear function as encountered
earlier in section 1.3, in definition \eref{stau}.)

The following proposition summarizes our findings, and proves (the case
$p=1$ is not conclusively proved by the variational equations above)
that we have indeed found the minimizer of
${\Phi}^{^{\SUR}}_{\mathbf{w},p}(f ; a)$:
\begin{proposition}
\label{prop-2-1} Suppose the operator $K$ maps a Hilbert space $\cH$
to another Hilbert space $\cH'$, with $\|\K K\| < 1$,
and suppose $g$ is an element of $\cH'$.
Let $(\vpg)_{\gamma \in \Gamma}$
be an orthonormal basis for $\cH$, and
let $\mathbf{w}=(\ag)_{\g \in \Gamma}$ be a sequence
of strictly positive numbers. Pick
arbitrary $p \geq 1$ and $a \in \cH$. Define the
functional ${\Phi}^{^{\SUR}}_{\mathbf{w},p}(f ; a)$ on $\cH$ by
$$
{\Phi}^{^{\SUR}}_{\mathbf{w},p}(f ; a)=\|Kf-g\|^2 + \sum_{\g \in
\Gamma} \ag |\fg|^p
+\|f-a\|^2-\|K(f-a)\|^2 ~.
$$
Then ${\Phi}^{^{\SUR}}_{\mathbf{w},p}(f ; a)$ has a unique minimizer
in $\cH$. \\
This minimizer
is given by $f=\S_{\mathbf{w},p}\left(a +\K (g-Ka) \right)$, where
the operators $\S_{\mathbf{w},p}$ are defined by
\begin{equation}
\label{def-SS}
\S_{\mathbf{w},p}(h)= \sum_{\gamma} S_{\ag,p}(h_{\gamma}) \vpg ~~,
\end{equation}
with the functions $S_{w,p}$ from $\R$ to itself given by
{\rm (\ref{S-pneq1}, \ref{S-peq1})}. For all $h \in \cH$, one has
$$
{\Phi}^{^{\SUR}}_{\mathbf{w},p}(f+h ; a)
\geq {\Phi}^{^{\SUR}}_{\mathbf{w},p}(f ;a)+\|h\|^2~.
$$
\end{proposition}
{\em Proof:}
The cases $p>1$ and $p=1$ should be treated slightly differently. We discuss
here only the case $p=1$; the simpler case $p>1$ is left to the reader.

Take $f'=f+h$, where $f$ is as defined in the Proposition, and $
h \in \cH$
is arbitrary. Then
$$
{\Phi}^{^{\SUR}}_{\mathbf{w},1}(f+h ; a) =
{\Phi}^{^{\SUR}}_{\mathbf{w},1}(f ; a)+ 2 \left<h,f-a -\K (g
-Ka) \right>
+\sum_{\g \in \Gamma} \ag (|\fg +h_{\g}| -|\fg|) + \|h\|^2  ~.
$$
Define now $\Gamma_{_{\!0}}=\{\g \in \Gamma; \fg=0\}$, and
$\Gamma_{_{\!1}}=\Gamma \setminus
\Gamma_{_{\!0}}$.
Substituting the explicit expression \eref{solcomp-peq1} for the $\fg$, we
have then
\begin{eqnarray*}
{\Phi}^{^{\SUR}}_{\mathbf{w},1}(f+h ; a)-
{\Phi}^{^{\SUR}}_{\mathbf{w},1}(f; a) &= &\|h\|^2
+ \sum_{\g \in \Gamma_{_{\!0}}}
\left[ \ag |h_{\g}| - 2 h_{\g} (a_{\g} +[ \K (g -K a)]_{\g} ) \right] \\
&&~~~~~~~~~~~~
+\sum_{\g \in \Gamma_{_{\!1}}} \left( \ag |\fg + h_{\g}| - \ag |\fg| + h_{\g}
[-\ag  ~ \mbox{sign}(\fg) ]\right) ~.
\end{eqnarray*}
For $\g \in \Gamma_{_{\!0}}$, $~2|a_{\g} +[ \K (g -K a)]_{\g}| \leq
\ag$, so
that
$\ag |h_{\g}| - 2 h_{\g}\, (a_{\g} +[ \K (g -K a)]_{\g}) \geq  0$.\\
If $\g \in \Gamma_{_{\!1}}$, we distinguish two cases, according to the sign
of $\fg$\, . If $\fg >0$, then\\
$\ag |\fg + h_{\g}| - \ag |\fg| + h_{\g}
[-\ag  ~ \mbox{sign}(\fg) ]= \ag [ |\fg +h_{\g}| - (\fg + h_{\g}) ] \geq 0$.
If $\fg <0$, then $\ag |\fg + h_{\g}| - \ag |\fg| + h_{\g}
[-\ag  ~ \mbox{sign}(\fg) ]= \ag [ |\fg +h_{\g}| + (\fg + h_{\g}) ] \geq 0$.\\
It follows that ${\Phi}^{^{\SUR}}_{\mathbf{w},1}(f+h;a)-
{\Phi}^{^{\SUR}}_{\mathbf{w},1}(f;a) \geq \|h\|^2 $, which
proves the Proposition.
 \hfill \QED

\bigskip

For later reference it is useful to point out that
\begin{lemma}
\label{SS-non-exp}
The operators $\S_{\mathbf{w},p}$ are non-expansive, i.e.
$\forall v, ~ v' \in \cH$, $\| \S_{\mathbf{w},p} v - \S_{\mathbf{w},p} v'
\| \leq
\|v-v'\|~.$
\end{lemma}
{\em Proof:}
As shown by \eref{def-SS},
$$
\|\S_{\mathbf{w},p} v - \S_{\mathbf{w},p} v' \|^2 =
\sum_{\g \in \Gamma} |S_{\ag,p}( v_{\g}) - S_{\ag,p}( v'_{\g})|^2 ~,
$$
which means that it suffices to show that, $\forall x , x' \in \R$, and all
$w\geq0$,
\begin{equation}
\label{S-non-exp}
| S_{w,p}(x)-S_{w,p}(x')| \leq |x-x'|~.
\end{equation}
If $p >1$, then $S_{w,p}$ is the inverse of the function $F_{w,p}$; since
$F_{w,p}$ is differentiable with derivative uniformly bounded below by 1,
\eref{S-non-exp} follows immediately in this case.\\
If $p=1$, then $S_{w,1}$ is not differentiable in $x=w/2$ or $x=-w/2$, and
another
argument must be used. For the sake of definiteness, let us assume $x \geq x'$.
We will just check all the possible cases.
If $x$ and $x'$ have the same sign and $|x|,~|x'|\geq w/2$, then
$| S_{w,p}(x)-S_{w,p}(x')| =|x-x'|$. If $x'\leq -w/2$ and $x \geq w/2$, then
$| S_{w,p}(x)-S_{w,p}(x')| = x +|x'|-w < |x-x'|$. If
$x \geq w/2$ and $|x'| < w/2$, then
$| S_{w,p}(x)-S_{w,p}(x')| =x-w/2 < |x-x'|$. A symmetric argument applies
to the
case $|x|<w/2$ and $x' \leq -w/2$. Finally, if both $|x|$ and $|x'|$ are
less than
$w/2$, we have
$| S_{w,p}(x)-S_{w,p}(x')|=0 \leq |x-x'|$. This establishes \eref{S-non-exp}
in all cases. \hfill \QED

\bigskip

Having found the minimizer of a generic
${\Phi}^{^{\SUR}}_{\mathbf{w},p}(f ; a)$, we can
apply this to the
iteration \eref{iter}, leading to

\begin{corollary}
\label{cor-2-2}
Let $\cH$, $\cH'$, $K$, $g$, $\mathbf{w}$ and $(\vpg)_{\g \in \Gamma}$ be
as in Proposition {\rm \ref{prop-2-1}}. Pick $f^0$ in $\cH$, and define
the functions $f^n$ recursively by the algorithm {\rm \eref{iter}}.
Then
\begin{equation}
\label{f-n}
f^n= \S_{\mathbf{w},p}\left(f^{n-1}+ \K (g-Kf^{n-1}) \right) ~.
\end{equation}
\end{corollary}
{\em Proof:} this follows immediately from Proposition \ref{prop-2-1}.
$~~~~~~~~~~~$ \hfill \QED

\begin{remark}
\label{op-D}
{\rm In the argument above, we used essentially only two ingredients: the
(strict)
convexity of $\|f-a\|^2-\|K(f-a)\|^2$, and the presence of the negative
$-\|Kf\|^2$
term in this expression, canceling the $\|Kf\|^2$ in the original
functional. We can use this observation to present a slight
generalization, in which the identity operator used to upper bound $ \K K
$ is replaced by a more general operator $D$ that is diagonal in the
$\vpg$--basis,}
$$
D\, \vpg = d_{\gamma} \vpg~,
$$
{\rm and that still gives a strict upper bound for $ \K K$, i.e. satisfies}
$$
D \geq K^*K + \eta I\ \  \mbox{\rm for some } \eta >0\ .
$$
{\rm In this case, the whole construction still carries through, with slight
modifications; the successive $f^n$ are now given by }
\begin{equation*}
\fg ^n = S_{\ag/d_{\gamma},p}\left(\fg ^{n-1} +
\frac {[\K (g-Kf^{n-1})]_{\gamma}}{d_{\gamma}} \right) ~~.
\end{equation*}
{\rm Introducing the notation $\mathbf{w/d}$ for the sequence
$(\ag/d_{\gamma})_{\gamma}$, we can rewrite this as }
\begin{equation*}
f^n = \S_{\mathbf{w/d},p}\left( f^{n-1} +D^{-1}[\K (g-Kf^{n-1})] \right) ~.
\end{equation*}
{\rm
For the sake of simplicity of notation,
we shall restrict ourselves
to the case $D = \mbox{\I}$.}
\end{remark}

\begin{remark}
\label{2-5}
{\rm If we deal with complex rather than real functions, and the $f_{\gamma}$,
$(K^*g)_{\gamma}, \cdots$ are complex quantities, then the derivation
of the minimizer of ${\Phi}^{^{\SUR}}_{\mathbf{w},1}(f; a)$
has to be adapted somewhat. Writing
$f_{\gamma}= r_{\gamma} e^{i \theta_{\gamma}}$, with
$ r_{\gamma} \geq 0$, $\theta_{\gamma} \in [0,2 \pi )$,
and likewise $(a +K^*g -K^*Ka)_{\gamma} = R_{\gamma} e^{i \Theta_{\gamma}}$,
we find, instead of \eref{sur},
$$
{\Phi}^{^{\SUR}}_{\mathbf{w},p}(f; a)
= \sum_{\gamma} [ r_{\gamma}^2+ w_{\gamma}r_{\gamma}^p
- 2r_{\gamma}R_{\gamma}\cos(\theta_{\gamma}-\Theta_{\gamma})]
+\|g\|^2+\|a\|^2-\|Ka\|^2~.
$$
Minimizing over $r_{\gamma} \in [0,\infty)$ and
$\theta_{\gamma} \in [0,2 \pi)$ leads to
$\theta_{\gamma}=\Theta_{\gamma}$ and
$r_{\gamma}=S_{w_{\gamma},p}(R_{\gamma})$.
If we extend the definition
of $S_{\mu,p}$ to complex arguments by setting
$S_{\mu,p}(r e^{i \theta}) = S_{\mu,p}(r) e^{i \theta}$,
then this still leads to
$\fg= S_{w_{\gamma},p}\left(a_{\gamma} +[K^*(g -Ka)]_{\gamma} \right)$,
as in (\ref{solcomp-pneq1}, \ref{solcomp-peq1}).
The arguments of the different proofs still
hold for this complex version, after minor and straightforward modifications.}
\end{remark}

\section{Convergence of the iterative algorithm}

In this section we discuss the convergence of the sequence $(f^n)_{n \in \N}$
defined by {\rm \eref{f-n}}. The main result of this section is the
following theorem:
\begin{theorem}
\label{th-3-1}
Let $K$ be a bounded linear operator from
$\cH$ to $\cH'$, with norm strictly bounded by $1$. Take $p \in [1,2]$, and
let $\S_{\mathbf{w},p}$ be the shrinkage operator defined by {\rm
\eref{def-SS}},
where the sequence $\mathbf{w}=(\ag)_{\g \in \Gamma}$ is uniformly
bounded
below away from zero, i.e. there exists a constant $c>0$ such that $\forall \g
\in
\Gamma:$
$\ag \geq c$.
Then the sequence of iterates
\begin{equation*}
f^n=\S_{\mathbf{w},p}\left( f^{n-1} + \K (g- Kf^{n-1})\right)\ ,\quad
n=1,2,\dots\;,
\end{equation*}
with $f^0$ arbitrarily chosen in $\cH$, converges strongly to a minimizer
of the functional
\begin{equation*}
\Phi_{\mathbf{w},p}(f) = \| Kf-g\|^2 +  \Vvert f\Vvert_{\mathbf{w},p}^p\ ,
\end{equation*}
where $\Vvert f\Vvert_{\mathbf{w},p}$ denotes the norm
\begin{equation}
\label{triple-norm}
\Vvert f\Vvert_{\mathbf{w},p} = \left[ \sum_{\g \in \Gamma} \ag
|\left<f,\vpg \right>|^p
\right]^{1/p} ~,~ 1 \leq p \leq 2~.
\end{equation}
If either $p>1$ or {\rm{N}}$(K)=\{0\}$, then the minimizer $f^\star$ of
$\Phi_{\mathbf{w},p}$ is unique, and every sequence of iterates
$f^n$ converges strongly to $f^\star$, regardless of the choice of $f^0$.
\end{theorem}

By ``strong convergence'' we mean convergence in the norm of $\cH$, as
opposed to
weak convergence.
This theorem will be proved in several stages. To start, we prove weak
convergence,
and we establish that the weak limit is indeed a minimizer of
$\Phi_{\mathbf{w},p}$.
Next, we prove that the convergence holds in norm, and not only in the weak
topology.
To lighten our formulas, we introduce the shorthand notation
$$
\T f = \S_{\mathbf{w},p}\left( f + \K (g-Kf)\right)~;
$$
with this new notation we have $f^n= \T^n f^0$.

\subsection{Weak convergence of the $f^n$}

To prove weak convergence of the $f^n=\T^n f^0$, we apply the following
theorem, due to
Opial \cite{Opi67}:
\begin{theorem}
\label{thm_opial}
Let
the mapping $\A $ from $\cH$ to $\cH$ satisfy the following
conditions:
\begin{enumerate}
\item[{\rm (i)}] $\A $ is non-expansive: $\forall v, v' \in \cH$,
$\| \A  v -  \A  v'\| \leq \| v - v'\|$,
\item[{\rm (ii)}] $\A $ is asymptotically regular: $\forall v \in \cH$,
$\| \A ^{n+1}v -\A ^n v\|
\xrightarrow[n \to \infty ]{~}  0$ ,
\item[{\rm (iii)}] the set ${\cal F}$ of the fixed points of $\A $ in
${\cH}$ is
not empty.
\end{enumerate}
Then, $\forall v \in \cH$, the sequence $(\A ^n v)_{n \in \N}$
converges weakly to a fixed point in ${\cal F}$.
\end{theorem}

Opial's original proof can be simplified;
we provide the simplified proof (still mainly
following Opial's approach) in Appendix \ref{Opial}.
(The theorem is slightly more general than what is stated
in Theorem \ref{thm_opial} in that the mapping $\A $ need not be
defined on all of space; it suffices that it map a closed convex subset of
$\cH$ to itself -- see Appendix \ref{Opial}.
\cite{Opi67} also contains additional refinements,
which we shall not need here.)
One of the Lemmas stated and proved in the Appendix
will be invoked in its own right, further below in this section;
for the reader's convenience, we state
it here in full as well:

\begin{lemma}
\label{lem-3-2}
Suppose the mapping $\A$ from $\cH$
to $\cH$ satisfies the conditions {\rm (i)} and {\rm (ii)}
in Theorem {\rm\ref{thm_opial}}. Then, if
a subsequence of $(\A ^n v)_{n\in \mathbb{N}}$
converges weakly in $\cH$, then its limit is a fixed point of $\A$.
\end{lemma}

In order to apply Opial's Theorem to our nonlinear operator
$\T$, we need to verify that it satisfies the three conditions in Theorem
\ref{thm_opial}. We do this in the following series of lemmas. We first have

\begin{lemma}
\label{nonexp}
The mapping $\T$ is non-expansive, i. e. $\forall v, v' \in \cH$
\begin{equation*}
\|\T  v -  \T  v' \| \leq \| v - v' \| \ .
\end{equation*}
\end{lemma}
{\em Proof:}
It follows from Lemma \ref{SS-non-exp}
that the shrinkage operator ${\S_{\mathbf{w},p}}$ is
non-expansive. Hence we have
\begin{eqnarray*}
\|{\T} v - {\T} v' \| &\leq& \| (I-K^*K) v - (I-K^*K)
v' \|\\
&\leq& \| I-K^*K  \|\ \| v - v' \|  \leq \| v - v' \|
\end{eqnarray*}
because we assumed $\| K \| < 1$.
\hfill\QED\bigskip

This verifies that $\T$ satisfies the first condition (i) in Theorem
\ref{thm_opial}.
To verify the second condition, we first prove some auxiliary lemmas.
\begin{lemma}
\label{cost2}
Both $\left({\Phi}_{\mathbf{w},p}(f^n)\right)_{n \in \N}$ and
$\left({\Phi}^{^{\SUR}}_{\mathbf{w},p}(f^{n+1} ; f^n)\right)_{n
\in
\N}$ are  non-increasing sequences.
\end{lemma}
{\em Proof:} For the sake of convenience, we introduce the operator
$L = \sqrt{I -\K K }$, so that $\|h\|^2-\|Kh\|^2= \|Lh\|^2$. Because
$f^{n+1}$ is the
minimizer of the functional
${\Phi}^{^{\SUR}}_{\mathbf{w},p}(f ; f^n)$ and therefore
\begin{equation*}
\Phi_{\mathbf{w},p}(f^{n+1})+\| L(f^{n+1}-f^n)\|^2 =
{\Phi}^{^{\SUR}}_{\mathbf{w},p}(f^{n+1} ; f^n)
\leq
{\Phi}^{^{\SUR}}_{\mathbf{w},p}(f^n ; f^n)=\Phi_{\mathbf{w},p}(f^n)\
,
\end{equation*}
we obtain
\begin{equation*}
\Phi_{\mathbf{w},p}(f^{n+1})\leq \Phi_{\mathbf{w},p}(f^n) \ .
\end{equation*}
On the other hand
\begin{equation*}
{\Phi}^{^{\SUR}}_{\mathbf{w},p}(f^{n+2} ; f^{n+1})\leq
\Phi_{\mathbf{w},p}(f^{n+1})
\leq \Phi_{\mathbf{w},p}(f^{n+1})+
\|L(f^{n+1}-f^n)\|^2={\Phi}^{^{\SUR}}_{\mathbf{w},p}(f^{n+1} ; f^n)\
.
\end{equation*}
\hfill\QED\bigskip

\begin{lemma}
\label{unifbddness}
Suppose the sequence $\mathbf{w}=(\ag)_{\g \in \Gamma}$ is uniformly bounded
below by a strictly positive number.
Then the $\|f^n\|$ are bounded uniformly in $n$.
\end{lemma}
{\em Proof:} Since $\ag \geq c$, uniformly in $\g$, for some $c>0$, we have
\begin{equation*}
\Vvert f^n\Vvert_{\mathbf{w},p}^p
\leq  \Phi_{\mathbf{w},p}(f^n) \leq
\Phi_{\mathbf{w},p}(f^0)~,
\end{equation*}
by Lemma \ref{cost2}. Hence the $f^n$ are bounded uniformly
in the $\Vvert ~\Vvert_{\mathbf{w},p}$-norm.
Since
\begin{equation}
\label{bdL2Ban}
\| f \|^2 \leq c^{-2/p} \ {\mathop{\rm max}_{\g \in \Gamma}}
[\ag^{(2-p)/p} |f_{\g}|^{2-p}]\ \Vvert f \Vvert_{\mathbf{w},p}^p
\leq c^{-2/p} \Vvert f\Vvert_{\mathbf{w},p}^{2-p}\ \Vvert f
\Vvert_{\mathbf{w},p}^p = c^{-2/p} \Vvert f
\Vvert_{\mathbf{w},p}^2\ ,
\end{equation}
we also have a uniform bound on the $\|f^n\|$.
\hfill\QED\bigskip

\begin{lemma}
\label{series}
The
series
$\sum_{n=0}^\infty \| f^{n+1}-f^n\|^2 $ is convergent.
\end{lemma}
{\em Proof:} This is a consequence of the strict positive-definiteness of $L$,
which holds because $\|K\| <1$. We have, for any $N \in \N$,
\begin{equation*}
\sum_{n=0}^N \| f^{n+1}-f^n\|^2 \leq \frac{1}{A} \sum_{n=0}^N
\| L(f^{n+1}-f^n)\|^2
\end{equation*}
where $A$ is a strictly positive lower bound for the spectrum of $L^*L$.
By Lemma \ref{cost2},
\begin{equation*}
\sum_{n=0}^{N}
\| L(f^{n+1}-f^n)\|^2 \leq
\sum_{n=0}^N [\Phi_{\mathbf{w},p}(f^n)-\Phi_{\mathbf{w},p}(f^{n+1})]
= \Phi_{\mathbf{w},p}(f^{0})-\Phi_{\mathbf{w},p}(f^{N+1}) \leq
\Phi_{\mathbf{w},p}(f^{0})~,
\end{equation*}
where we have used that
$(\Phi_{\mathbf{w},p}(f^n))_{n \in \N}$ is a non-increasing sequence. \\
It follows
that $\sum_{n=0}^N \| f^{n+1}-f^n\|^2 $ is bounded uniformly in $N$, so that
the infinite series converges.\hfill \QED

\bigskip

As an immediate consequence, we have that
\begin{lemma}
\label{asyreg}
The mapping ${\T}$ is asymptotically regular, i.e.
\begin{equation*}
\|{\T}^{n+1}f^0 -{\T}^n f^0 \| =
\| f^{n+1} - f^n \| \to 0 \quad {\rm for} \quad n \to \infty\ .
\end{equation*}
\end{lemma}
We can now establish the following
\begin{proposition}
The sequence $f^n={\T}^nf^0$,
$n=1, 2, \cdots$ converges weakly, and its limit is a fixed point for $\T$.
\end{proposition}
{\em Proof:} Since, by Lemma \ref{unifbddness}, the $f^n={\T}^n f^0$
are uniformly bounded in $n$, it follows from the Banach-Alaoglu
theorem that they have a weak accumulation point. By Lemma \ref{lem-3-2},
this weak accumulation point is a fixed point for $\T$.
It follows that the set of fixed points of $\T$ is not empty.
Since $\T$ is also non-expansive (by Lemma \ref{nonexp}) and
asymptotically regular (by Lemma \ref{asyreg}),
we can apply Opial's theorem (Theorem \ref{th-3-1} above), and the
conclusion of the Proposition follows.
\hfill\QED

\bigskip

By the following proposition this fixed point is also a minimizer
for the functional $\Phi_{\mathbf{w},p}$.
\begin{proposition}
\label{fix-min}
A fixed point for ${\T}$ is a minimizer for the functional
$\Phi_{\mathbf{w},p}$.
\end{proposition}
{\em Proof:}
If $f^\star = {\T} f^\star$, then by Proposition
\ref{prop-2-1}, we know that $f^\star$ is a minimizer for the surrogate
functional ${\Phi}^{^{\SUR}}_{\mathbf{w},p}(f ; f^\star)$, and
that, $\forall h \in \cH$,
\begin{equation*}
{\Phi}^{^{\SUR}}_{\mathbf{w},p}(f^\star + h ; f^\star) \geq
{\Phi}^{^{\SUR}}_{\mathbf{w},p}(f^\star ;f^\star) + \| h \|^2
~.
\end{equation*}
Observing that ${\Phi}^{^{\SUR}}_{\mathbf{w},p}(f^\star ; f^\star) =
\Phi_{\mathbf{w},p}(f^\star)$, and
\begin{equation*}
{\Phi}^{^{\SUR}}_{\mathbf{w},p}(f^\star + h ; f^\star) =
\Phi_{\mathbf{w},p}(f^\star + h) + \| h \|^2 - \| Kh\|^2 \ ,
\end{equation*}
we conclude that, $\forall h \in \cH$,
$ \Phi_{\mathbf{w},p}(f^\star + h) \geq \Phi_{\mathbf{w},p}(f^\star) + \|
Kh\|^2$, which shows that $f^\star$ is a minimizer for $\Phi(f)$.
\hfill\QED

\bigskip

The following proposition summarizes this subsection.
\begin{proposition}
\label{prop-wk-conv}
{\rm (Weak convergence)} Make the same assumptions as in the statement of
{\rm Theorem
\ref{th-3-1}}. Then, for
any choice of the initial $f^0$, the sequence $f^n={\T}^n f^0, \
n=1,2,\cdots$ converges weakly to a minimizer for $\Phi_{\mathbf{w},p}$.
If either {\rm N}$ (K)=\{0\}$ or $p>1$, then $\Phi_{\mathbf{w},p}$ has a unique
minimizer $f^\star$, and all the sequences $(f^n)_{n \in \N}$ converge
weakly to $f^\star$, regardless of the choice of $f^0$.
\end{proposition}
{\em Proof:}
The only thing that hasn't been proved yet above is the
uniqueness of the minimizer if $\mN (K)=\{0\}$ or $p>1$. This uniqueness
follows from the observation that $\Vvert f \Vvert_{\mathbf{w},p}$ is strictly
convex in $f$ if $p>1$, and that $\|Kf-g\|^2$ is strictly convex in
$f$ if $\mN (K)=\{0\}$. In both these cases $\Phi_{\mathbf{w},p}$ is thus
strictly convex, so that it has a unique minimizer. \hfill \QED

\bigskip

\begin{remark}
{\rm If one has the additional prior information that the object lies in
some closed convex subset ${\cal C}$ of the Hilbert space $\cH$, then the
iterative procedure can be adapted to take this into account, by replacing
the shrinkage operator
${\S}$  by ${\mathbf P}_{\!\cal C} {\S}$,
where ${\mathbf P}_{\!\cal C}$ is the projector on ${\cal C}$. For example,
if $\cH=L^2$, then ${\cal C}$ could be the cone of functions that are positive
almost everywhere. The  results in this subsection can be extended to this
case;
a more general version of Theorem \ref{thm_opial}
can be applied, in which $\A $ need not be defined
on all of $\cH$, but only on ${\cal C} \subset \cH$; see Appendix \ref{Opial}.
We would however need to either use other tools
to ensure, or assume outright
that the set of fixed points of $\T={\mathbf P}_{\!\cal C} {\S}$ is not empty
(see also  \cite{Eic92})}.
\end{remark}

\bigskip

\begin{remark}
{\rm If $\Phi_{\mathbf{w},p}$ is strictly convex, then one can prove the weak
convergence more directly, as follows. By the boundedness of the $f^n$
(Lemma \ref{unifbddness}), we must have a weakly convergent subsequence
$(f^{n_k})_{k \in \N}$. By Lemma \ref{asyreg}, the sequence
$(f^{n_k+1})_{k \in \N}$ must then also be weakly convergent, with the same
weak limit $\wf$. It then follows from the equation
$$
f^{n_k+1}_\g= S_{\ag,p}\left(f^{n_k}_\g +[\K (g-K f^{n_k})]_\g \right)~,
$$
together with $\lim_{k \to \infty}f^{n_k}_\g = \lim_{k \to \infty}f^{n_k+1}_\g
=\wf _\g$, that $\wf$ must be the fixed point $f^\star$ of T. Since this
holds for
any weak accumulation point of $(f^n)_{n \in \N}$, the weak convergence
of $(f^n)_{n \in \N}$ to $f^\star$ follows. }
\end{remark}

\bigskip

\begin{remark}
{\rm The proof of Lemma \ref{unifbddness} is the only place, so far, where we
have explicitly used $p \leq 2$. If it were possible to establish a uniform
bound on the $\|f^n\|$ by some other means (e.g. by showing that the
$\|\T^n f^0\|$
are bounded uniformly in $n$), then we could dispense with the restriction
$p \leq 2$, and Proposition \ref{prop-wk-conv} would hold for all $p \geq 1$. }
\end{remark}

\subsection{Strong convergence of the $f^n$}

In this subsection we shall prove that the convergence of the successive
iterates
$\{f^n\}$ holds not only in the weak topology, but also in the Hilbert
space norm.
Again, we break up the proof into several lemmas. For the sake of convenience,
we introduce the following notations
\begin{eqnarray}
f^\star&=& \mbox{{\em w}\! --\!}\lim_{n \to \infty} f^n \nonumber \\
u^n &=& f^n - f^\star \nonumber \\ 
h\ &=& f^\star + \K (g-Kf^\star)\ . \label{redef2}
\end{eqnarray}
Here and below, we use the notation {\em w}$\,$--$\lim$ as a shorthand
for {\em weak limit}.
\begin{lemma}
\label{Ku}
$\| Ku^n \| \to 0$ for $n \to \infty$\ .
\end{lemma}
{\em Proof:}
Since
\begin{equation*}
u^{n+1} - u^n = {\S}_{\mathbf{w},p}\left( h+(I-K^*K)u^n\right) -
{\S}_{\mathbf{w},p}(h) - u^n
\end{equation*}
and $\|u^{n+1} - u^n \| = \| f^{n+1} - f^n \| \to 0\ {\rm for}\ n
\to \infty$ by Lemma \ref{asyreg}, we have
\begin{equation}
\|\ {\S}_{\mathbf{w},p}\left( h+(I-K^*K)u^n\right) - {\S}_{\mathbf{w},p}(h) -
u^n \| \to 0 \  {\rm for}\ n \to \infty ~,
\label{cvaux}
\end{equation}
and hence also
\begin{equation}
 \max\left(0,\| u^n \|- \|{\S}_{\mathbf{w},p}\left( h+(I-K^*K)u^n\right) -
{\S}_{\mathbf{w},p}(h)\|\ \right) \to 0 \ {\rm for}\  n \to \infty\ .
\label{cvtozero}
\end{equation}
Since ${\S}_{\mathbf{w},p}$ is non-expansive
(Lemma \ref{SS-non-exp}), we have
\begin{equation*}
\| \ {\S_{\mathbf{w},p}}\left( h+(I-K^*K)u^n\right) - {\S}_{\mathbf{w},p}(h) \|
\ \leq
\| (I-K^*K)u^n \| \leq \| u^n \|~;
\end{equation*}
therefore the ``max'' in  (\ref{cvtozero}) can be dropped, and it follows that
\begin{equation}
\| u^n \| - \| (I-K^*K)u^n \|  \to 0 \ {\rm for}\  n \to \infty \ .
\label{cv2}
\end{equation}
Because
\begin{eqnarray*}
\| u^n \| + \| (I-K^*K)u^n \| &\leq& 2\| u^n\| = 2\|f^n-f^\star\|\\
&\leq& 2(\| f^\star\|+ {\mathop{\rm sup}_{k}} \| f^k\|) = C
\end{eqnarray*}
where $C$ is a finite constant (by Lemma \ref{unifbddness}), we obtain
\begin{equation*}
0 \leq \| u^n \|^2 - \| (I-K^*K)u^n \|^2 \leq
C (\| u^n \| - \| (I-K^*K)u^n \|)\ ,
\end{equation*}
which tends to zero by (\ref{cv2}).
The inequality
\begin{equation*}
\| u^n \|^2 - \| (I-K^*K)u^n \|^2 =
2\| Ku^n\|^2-\|\K Ku^n\|^2 \geq \| Ku^n\|^2
\end{equation*}
then implies that $\| Ku^n\|^2 \to 0 \ {\rm for}\  n \to \infty\ $.
\hfill\QED
\begin{remark}
\label{ifKcomp}
{\rm Note that if $K$ is a compact operator, the weak convergence
to $0$ of the $u_n$ automatically implies that $\|K u_n\|$ tends
to $0$ as $n$ tends to $\infty$, so that we don't need
Lemma \ref{Ku} in this case.}
\end{remark}

\bigskip

If $K$ had a bounded inverse, we could conclude from $\|K u_n\| \to 0$ that
$\| u_n\|
\to 0
\ {\rm for}\  n \to \infty\ $. If this is not the case, however,
and thus for all
ill-posed linear inverse problems, we need some extra work to show the norm
convergence of $f^n$ to $f^\star$.
\begin{lemma}
For $h$ given by {\rm (\ref{redef2})},
$\| {\S}_{\mathbf{w},p}(h+u^n) - {\S}_{\mathbf{w},p}(h) - u^n \| \to 0$ for
$n \to
\infty$.
\end{lemma}
{\em Proof:}
We have
\begin{eqnarray*}
\| {\S}_{\mathbf{w},p}(h+u^n) - {\S}_{\mathbf{w},p}(h) - u^n \|
&\leq& \| {\S}_{\mathbf{w},p}(h+u^n-K^*Ku^n) - {\S}_{\mathbf{w},p}(h) - u^n
\|\\ &&~~~~+
\|{\S}_{\mathbf{w},p}(h+u^n) - {\S}_{\mathbf{w},p}(h+u^n-K^*Ku^n) \| \\
&\leq& \| {\S}_{\mathbf{w},p}(h+u^n-K^*Ku^n) - {\S}_{\mathbf{w},p}(h) - u^n
\|\\ &&
~~~~+\|\K Ku^n \|~,
\end{eqnarray*}
where we used the non-expansivity of ${\S}_{\mathbf{w},p}$ (Lemma
\ref{SS-non-exp}).
The result follows since both terms in this last bound tend to zero for $n \to
\infty$ because of Lemma
\ref{Ku} and (\ref{cvaux}).
\hfill\QED

\bigskip

\begin{lemma}
\label{lm-3-16}
If for some $a \in \cH$, and some sequence $(v^n)_{n \in \N}$,
w--$\lim_{n \to \infty}v^n=0$  and
$\lim_{n \to \infty}\| { \S}_{\mathbf{w},p}(a+v^n) -
{\S}_{\mathbf{w},p}(a) - v^n \|=0$
then $\| v^n \| \to 0$ for $n \to \infty$.
\end{lemma}
{\em Proof:}
The argument of the proof is slightly different for the cases $p=1$
and $p>1$, and we treat the two cases separately. \\
We start with $p>1$.
Since the sequence $\{v^n\}$ is weakly convergent, it has to be bounded: there
is a constant $B$ such that $\forall n$, $\| v^n \| \leq B$, and
hence also $\forall n, \forall \g \in \Gamma$, $\vert v^n_\g \vert \leq B$.
Next, we define the set $\Gamma_{_{\!0}}
= \{ \g \in \Gamma; |a_{\g}| \geq B\}$; since $a \in \cH$, this is a finite
set. We then have $\forall \g \in \Gamma_{_{\!1}}=\Gamma \setminus
\Gamma_{_{\!0}}$,
that $|a_{\g}|$ and $|a_{\g}+v^n_{\g}|$ are bounded above by $2B$.
Recalling the definition of $S_{\ag,p}=\left(F_{\ag,p}\right)^{-1}$,
and observing that, because $p \leq 2$,
$F'_{\ag,p}(x)=1+\ag p(p-1)|x|^{p-2}/2 \geq
1+ \ag \, p(p-1) /[2(2B)^{2-p}]$ if $|x| \leq 2B$ , we have
\begin{eqnarray*}
|S_{\ag,p}(a_{\g}+v^n_{\g}) - S_{\ag,p}(a_{\g})|
&\leq &\left( {\mathop {\max}_{|x|\leq 2B}} |S'_{\ag,p}(x)| \right) |v^n_{\g}|
\\ & \leq & \left( 1+ \ag\, p(p-1) /[2 (2B)^{2-p}] \right)^{-1} |v^n_{\g}| \\
& \leq & \left( 1+ c \, p(p-1) / [2(2B)^{2-p}] \right)^{-1} |v^n_{\g}|~;
\end{eqnarray*}
in the second inequality,
we have used that $|S_{\ag,p}(x)| \leq |x|$, a consequence of the
non-expansivity
 of $S_{\ag,p}$ (see
Lemma
\ref{SS-non-exp}) to upper bound the derivative
$S'_{\ag,p}$ on the interval $[-2B,2B]$ by the inverse of the lower bound for
$F'_{\ag ,p}$ on the same interval;
in the last inequality we used the uniform lower bound on the $\ag$, i.e.
$\forall \g,
~ \ag \geq c >0$.
Rewriting $\left( 1+ c \, p(p-1) / [2(2B)^{2-p}] \right)^{-1}= C'<1$, we
have thus,
$\forall \g \in \Gamma_{_{\!1}}$, $C'|v^n_{\g}| \geq
|S_{\ag,p}(a_{\g}+v^n_{\g}) - S_{\ag,p}(a_{\g})|$, which implies
$$
\sum_{\g \in \Gamma_{_{\!1}}} |v^n_{\g}|^2 \leq
\frac{1}{(1-C')^2} \sum_{\g \in \Gamma_{_{\!1}}} |v^n_{\g}
-S_{\ag,p}(a_{\g}+v^n_{\g})
+ S_{\ag,p}(a_{\g})|^2  \to 0 \mbox{ as } n \to \infty ~.
$$
On the other hand, since $\Gamma_{_{\!0}}$ is a finite set, and the $v^n$
tend to
zero weakly as $n$ tends to $\infty$, we also have
$$
\sum_{\g \in \Gamma_{_{\!0}}} |v^n_{\g}|^2 \to \infty \mbox{ as } n \to
\infty ~.
$$
This proves the proposition for the case $p>1$. \\
For $p=1$,
we define  a finite set $\Gamma_{_{\!0}} \subset \Gamma$ so
that $\sum_{ \g \in \Gamma \setminus \Gamma_{_{\!0}}}
|a_{\g}|^2 \leq (c/4 )^2$,
where $c$ is again the uniform lower bound on the $\ag$.
Because this is a finite set, the weak convergence of the $v^n$
implies that $\sum_{\g \in \Gamma_{_{\!0}}} |v^n_{\g}|^2
\xrightarrow[n \to \infty]{~}  0$,
so that we can concentrate on
$\sum_{\g \in \Gamma \setminus \Gamma_{_{\!0}}} |v^n_{\g}|^2$ only. \\
For each $n$, we split $\Gamma_{_{\!1}}=\Gamma \setminus \Gamma_{_{\!0}}$ into
two subsets:
$\Gamma_{_{\!1,n}} = \{\g \in \Gamma_{_{\!1}};
|v^n_{\g}+a_{\g}| <\ag/2\}$ and $\widetilde{\Gamma}_{_{\!1,n}}=
\Gamma_{_{\!1}} \setminus \Gamma_{_{\!1,n}}$. If $\g \in \Gamma_{_{\!1,n}}$,
then $S_{\ag,1}(a_{\g}+v^n_{\g})=
S_{\ag,1}(a_{\g}) =0$ (since $|a_{\g}|\leq c/4 \leq \ag/2$),
so that $|v^n_{\g} -
S_{\ag,1}(a_{\g}+v^n_{\g}) + S_{\ag,1}(a_{\g})|=|v^n_{\g}|$.
It follows that
$$
\sum_{\g \in \Gamma_{_{\!1,n}}} |v^n_{\g}|^2
\leq \sum_{\g \in \Gamma} |v^n_{\g} -
S_{\ag,1}(a_{\g}+v^n_{\g}) + S_{\ag,1}(a_{\g})|^2 \to 0 \mbox{ as }
n \to \infty ~.
$$
It remains to prove only that
the remaining sum, $\sum_{\g \in \widetilde{\Gamma}_{_{\!1,n}}}
|v^n_{\g}|^2 $ also tends
to $0$ as $n \to \infty$.  \\
If $\g \in \Gamma_{_{\!1}}$ and $|v^n_{\g}+a_{\g}| \geq \ag/2$, then
$|v^n_{\g}|\geq |v^n_{\g}+a_{\g}| - |a_{\g}| \geq \ag/2 -c/4
\geq c/4 \geq |a_{\g}|$, so that $ v^n_{\g}+a_{\g}$ and
$v^n_{\g}$ have the same sign; it then follows that
\begin{eqnarray*}
&& |v^n_{\g} -
S_{\ag,1}(a_{\g}+v^n_{\g}) + S_{\ag,1}(a_{\g})|=
|v^n_{\g} -
S_{\ag,1}(a_{\g}+v^n_{\g})| \\
&&~~~~~~~~
=|v^n_{\g}- (a_{\g}+v^n_{\g})+ \frac{\ag}{2} \mbox{sign}(v^n_{\g})|
\geq  \frac{\ag}{2} -|a_{\g}| \geq \frac{c}{4} ~.
\end{eqnarray*}
This implies that
$$
\sum_{\g \in \widetilde{\Gamma}_{_{\!1,n}}}|v^n_{\g} -
S_{\ag,1}(a_{\g}+v^n_{\g}) + S_{\ag,1}(a_{\g})|^2
\geq \left(\frac{c}{4}\right)^2 \# \widetilde{\Gamma}_{_{\!1,n}} ~;
$$
since $\|v^n-\S_{\mathbf{w},1}(a+v^n)+\S_{\mathbf{w},1}(a)\|
\xrightarrow[n \to \infty]{~} 0$, we know on the other hand
that
$$
\sum_{\g \in \widetilde{\Gamma}_{_{\!1,n}}}|v^n_{\g} -
S_{\ag,1}(a_{\g}+v^n_{\g}) + S_{\ag,1}(a_{\g})|^2 <
\left(\frac{c}{4}\right)^2
$$
when $n$ exceeds some threshold, $N$, which implies that
$\widetilde{\Gamma}_{_{\!1,n}}$ is empty when $n > N$. Consequently
$\sum_{\g \in \widetilde{\Gamma}_{_{\!1,n}}}|v^n_{\g}|^2 =0$
for $n > N$. This completes the proof for the case $p=1$.\hfill \QED

\bigskip

Combining the Lemmas in this subsection with the results of the previous
subsection gives a complete proof of Theorem \ref{th-3-1} as stated at the
start of this section.

\section{Regularization properties and stability estimates}

In the preceding section, we devised
an iterative algorithm that
converges towards a minimizer of the functional
\begin{equation}
\Phi_{\mathbf{w},p}(f) = \| Kf-g\|^2 +  \Vvert f\Vvert^p_{\mathbf{w},p}~.
\label{phimu2}
\end{equation}
For simplicity, let us assume, until further notice,
that either $p>1$ or $\mN(K)=\{0\}$, so that there is a unique minimizer.

In this section, we shall discuss to what extent this minimizer is
acceptable as a {\em regularized solution} of the (possibly ill-posed)
inverse problem $Kf=g$. Of particular interest to us is the {\em stability}
of the estimate. For instance, if $\mN(K)=\{0\}$, we would like to know
 to what extent the proposed solution, in this
case the minimizer of $\Phi_{\mathbf{w},p}$, deviates from the ideal
solution $f_o$ if the data are a (small) perturbation of the image
$Kf_o$ of $f_o$. (If $\mN(K) \neq \{0\}$, then there exist other $f$ that
have the same image as $f_o$, and the algorithm might choose one of those
--  see below.) In this discussion both the ``size'' of the perturbation
and the weight of the penalty term in the variational functional, given
by the coefficients $(\ag)_{\g \in \Gamma}$, play a role. We argued earlier
that we need
$\mathbf{w} \neq \mathbf{0}$ in order to provide a meaningful estimate
if e.g. $K$ is a compact operator; on the other hand,
if $g = Kf_o$, then the presence of the penalty term will cause the
minimizer of $\Phi_{\mathbf{w},p}$ to be different from $f_o$. We therefore
need to strike a balance between the respective weights of the
perturbation $g - Kf_o$ and
the penalty term. Let us first define a framework in which we can make this
statement more precise.

Because we shall deal in this section with data functions $g$ that are not
fixed, we
adjust our notation for the variational functional to make the dependence
on $g$
explicit
where appropriate: with this more elaborate notation, the right hand side
of, for instance, \eref{phimu2} is now $\Phi_{\mathbf{w},p; g}(f )$.
(Because we work with one fixed
operator $K$,  the dependence of the functional on $K$ remains ``silent''.)
In order to make it possible to vary the weight of the penalty term in the
functional, we introduce an extra parameter $\mu$. We shall thus consider
the functional
\begin{equation}
\label{phimu}
\Phi_{\mu,\mathbf{w},p; g}(f)=
\| Kf-g\|^2 +\mu \Vvert f\Vvert^p_{\mathbf{w},p}~.
\end{equation}
Its minimizer will likewise depend on all these parameters. In its full
glory, we denote it by $f^{\star}_{\mu ,\mathbf{w},p;g}$; when  confusion
is impossible we  abbreviate this notation. In particular, since $\mathbf{w}$
and $p$ typically will not vary in the limit procedure that defines stability,
we may omit them in the heat of the discussion. Notice that the dependence on
$\mathbf{w}$ and $\mu$ arises only through the product $\mu \mathbf{w}$.

As mentioned above, if the ``error'' $e =g - Kf_o$ tends to zero, we would like
to see our estimate for the solution of the inverse problem tend to $f_o$;
since the minimizer of $\Phi_{\mu ,\mathbf{w},p; g}(f)$ differs from $f_o$ if
$\mu \ne 0$, this means
that we shall have to consider simultaneously a limit for $\mu \rightarrow 0$.
More precisely, we want to find a functional dependence of $\mu$
on the noise level $\epsilon$, $\mu=\mu(\epsilon)$
such that
\begin{equation}
\label{desired-res}
\mu(\epsilon)  \xrightarrow[\epsilon \rightarrow 0]{~}  0  ~~~
\mbox {and } ~~~ \sup_{\|g-Kf_o\|\leq \epsilon}
\|f^{\star}_{\mu(\epsilon) ,\mathbf{w},p;g}-f_o \|
\xrightarrow[\epsilon \rightarrow 0]{~}  0~.
\end{equation}
for each $f_o$ in a certain class of functions.
If we can achieve this, then the ill-posed inverse problem will be {\em
regularized}
(in norm or ``strongly'') by our iterative method,
and $f^\star_{\mu,\mathbf{w},p;g}$ will be
called a {\em regularized solution}.
One  also says in this case
that the minimization of the penalized least-squares functional
(\ref{phimu2}) provides us with a {\em regularizing algorithm} or
{\em regularization method}.

\subsection {A general regularization theorem}
If the $\ag$ tend to $\infty$, or more precisely, if
\begin{equation}
\label{compact-emb}
\forall C >0 ~: ~ \# \{ \g \in \Gamma ; \ag \leq C \} < \infty ~,
\end{equation}
then the embedding of $\mathcal{B}_{\mathbf{w},p}=
\{ f \in \cH;\sum_{\g \in \Gamma} \ag |f_\g|^p < \infty \}$ in $\cH$ is
compact. (This is because the identity operator from
$\mathcal{B}_{\mathbf{w},p}$ to
$\cH$ is then the norm--limit in $\mathcal{L}(\mathcal{B}_{\mathbf{w},p},
\cH)$,
as $C \to \infty$,
of the finite rank operators $P_C$ defined by
$P_C f=\sum_{\g \in \Gamma_C} \ag \left<f,\varphi_\g \right> \varphi_\g$,
where $\Gamma_C= \{ \g \in \Gamma ; \ag \leq C \}$.) In this case,
general compactness arguments can be used to show that \eref{desired-res}
can be achieved. (See also further below.)
We are, however, also interested in the general case, where
the $\ag$ need not grow unboundedly.
The following theorem proves that we can then nevertheless
choose the dependence $\mu(\epsilon)$
so that \eref{desired-res} holds:
\begin{theorem}
\label{regthm}
Assume that $K$ is a bounded operator from $\cH$ to $\cH'$ with $\|K\|<1$, that
$1 \leq p \leq 2$ and that the entries in the sequence
$\mathbf{w}=(\ag)_{\g \in \Gamma}$
are bounded below uniformly by a strictly positive number $c$.
Assume that either $p>1$ or $\mbox{{\rm N}}(K)=\{0\}$.
For any $g \in \cH'$
and any $\mu >0$, define
$f^{\star}_{\mu, \mathbf{w},p;g}$ to be the minimizer of $\Phi_{\mu,
\mathbf{w},p ; g}(f)$.
If $\mu=\mu(\epsilon)$ satisfies the requirements
\begin{equation}
\label{mu-req}
\lim_{\epsilon \rightarrow 0} \mu(\epsilon)=0 ~~~~~ \mbox{{\rm and}}
~~~~~ \lim_{\epsilon \rightarrow 0} \epsilon^2/\mu(\epsilon) =0 ~,
\end{equation}
then we have, for any $f_o \in \cH$,
$$
\lim_{\epsilon \rightarrow 0} \left[ \sup_{\|g-Kf_o\|\leq \epsilon}
\|f^{\star}_{\mu(\epsilon), \mathbf{w},p;g}-\fs \|\right] =0 ~,
$$
where $\fs$ is the unique element of minimum
$\Vvert ~ \Vvert _{\mathbf{w},p}$--norm in $\mathcal{S}
= \mbox{{\rm N}}(K)+f_o =
\{f; Kf = Kf_o\}$.
\end{theorem}

Note that under the conditions of Theorem \ref{regthm}, $\fs$ must indeed
be unique:
if $p>1$, then the $\Vvert ~ \Vvert _{\mathbf{w},p}$--norm is strictly
convex, so
that there is a unique minimizer for this norm in the hyperspace $\mN(K)+f_o$;
if $p=1$, our assumptions require $\mN(K)=\{0\}$. Note also that
if $\mN(K)=\{0\}$ (whether or not $p=1$), then necessarily $\fs = f_o$.

To prove Theorem \ref{regthm}, we will need the following two lemmas:

\begin{lemma}
\label{lm-4-1}
The functions $S_{w,p}$ from $\R$ to itself, defined by {\rm (\ref{S-pneq1},
\ref{S-peq1})} for $p>1$, $p=1$, respectively, satisfy
$$
|S_{w,p}(x)-x| \leq {\frac{wp}{2}}\ |x|^{p-1}~.
$$
\end{lemma}
{\em Proof:}
For $p=1$, the definition \eref{S-peq1} implies immediately that
$|x-S_{w,1}(x)|= \min(w/2,|x|) \leq w/2$, so that the proposition holds
for $x \neq 0$. For $x=0$, $S_{w,1}(x)=0$.

For $p>1$, $S_{w,p}= \left( F_{w,p} \right)^{-1}$, where
$F_{w,p}(y)=y+{\frac{wp}{2}}|y|^{p-1} \mbox{sign}(y)$ satisfies $|F_{w,p}(y)|
\geq |y|$, and $|F_{w,p}(y)-y| \leq {\frac{wp}{2}} |y|^{p-1}$.
It follows that $\, |S_{w,p}(x)| \leq |x|\,$, and
$\, |x-S_{w,p}(x)| $ $\leq {\frac{wp}{2}} \; |S_{w,p}(x)|^{p-1} $
$ \leq {\frac{wp}{2}}\; |x|^{p-1} ~$.
$~~~~~~~~~~~~~~~~~~$ \hfill \QED

\bigskip

\begin{lemma}
\label{lm-conv}
If the sequence of vectors $\left(v_k\right)_{ _{k \in \N}}$ converges weakly
in $\cH$ to $v$, and $\lim_{k \to \infty} \Vvert v_k \Vvert_{\mathbf{w},p}$
$ = \Vvert v \Vvert_{\mathbf{w},p}$,
then $\left(v_k\right)_{ _{k \in \N}}$ converges
to $v$ in the $\cH$--norm, i.e. $\lim_{k \to \infty} \|v-v_k\|=0~$.
\end{lemma}
{\em Proof:}
It is a standard result that if {\em w}$\,$--$\lim_{k \to \infty} v_k = v$,
and
$\lim_{k \to \infty} \|v_k\|=\|v\|$, then $\lim_{k \to \infty} \|v- v_k\|^2
= \lim_{k \to \infty} \left( \|v\|^2 + \|v_k\|^2 - 2 \left< v, v_k \right>
\right)
=  \|v\|^2 + \|v\|^2 - 2 \left< v, v \right>
 = 0$.
We thus need to prove only that $\lim_{k \to \infty} \|v_k\|=\|v\|$.

Since the $v_k$ converge weakly, they are uniformly bounded. It follows that
the $|v_{k,\g}| = |\left<v_k, \vpg \right>|$ are bounded uniformly in $k$
and $\g$
by some finite number $C$. Define $r=2/p$. Since, for $x, y > 0$,
$|x^r - y^r| \leq r |x-y| \max(x,y)^{r-1}$, it follows that
$ \left|~|v_{k,\g}|^2 -|v_{\g}|^2 \right|
\leq r \, C^{p(r-1)}~  \left| ~ |v_{k,\g}|^p -|v_{\g}|^p \right|~.$
Because the $\ag$ are uniformly bounded below by $c>0$, we obtain
$$
\left| \|v_k\|^2 -\|v\|^2 \right| \leq
\sum_{\g \in \Gamma} \left| |v_{k,\g}|^2 - |v_{\g}|^2 \right|
\leq \frac{2}{c \, p}\, C^{2-p} \sum_{\g \in \Gamma}
\ag \left| ~ |v_{k,\g}|^p -|v_{\g}|^p \right| ~,
$$
so that it suffices to prove that this last expression tends to $0$
as $k$ tends to $\infty$.
Define now $u_{k,\g}=\min \left( |v_{k,\g}|,|v_{\g}| \right)$. Clearly
$ \forall \g \in \Gamma ~:~ \lim_{k \to \infty} u_{k,\g}= |v_{\g}|~$; since
$\sum_{\g \in \Gamma}  \ag |v_{\g}|^p < \infty$, it follows by the dominated
convergence theorem that $\lim_{k \to \infty}
\sum_{\g \in \Gamma}  \ag u_{k,\g}^p =
\sum_{\g \in \Gamma}  \ag |v_{\g}|^p $. Since
$$
\sum_{\g \in \Gamma}
\ag \left| ~ |v_{k,\g}|^p -|v_{\g}|^p \right| =
\sum_{\g \in \Gamma}  \ag \left( |v_{\g}|^p + |v_{k,\g}|^p
- 2 u_{k,\g}^p \right) \xrightarrow[k \to \infty]{~} 0 ~,
$$
the Lemma follows. \hfill \QED

We are now ready to proceed to the

{\em Proof of }Theorem 4.1:
\newline
Let's assume that $\mu(\epsilon)$ satisfies the requirements \eref{mu-req}.
\newline
We first establish weak convergence. For this it is sufficient to prove
that if $(g_n)_{n \in \N}$ is a sequence in $\cH'$ such that
$\|g_n-Kf_o\| \leq \epsilon_n$, where
$(\epsilon_n)_{n \in \N}$ is a sequence of strictly positive numbers
that converges to zero
as $n \to \infty$, then {\em w}$\,$--$\lim_{n \to \infty}
f^{\star}_{\mu(\epsilon_n);g_n}= f^\dagger$, where $f^{\star}_{\mu;g}$
is the unique minimizer of $\Phi_{\mu ,\mathbf{w},p;g}(f)$
(As predicted, we have dropped here the
explicit indication of the dependence of $f^{\star}$ on $\mathbf{w}$ and
$p$; these
parameters will keep fixed values throughout this proof. We will take the
liberty
to drop them in our notation for $\Phi$ as well, when this is convenient.)
For the sake of convenience,
we abbreviate $\mu(\epsilon_n)$ as $\mu_n$. \\
Then the $f^{\star}_{\mu_n;g_n}$ are uniformly bounded in $\cH$
by the following argument:
\begin{eqnarray}
\|  f^\star_{\mu_n;g_n} \|^p &\leq &  \frac{1}{c}
\Vvert f^\star_{\mu_n;g_n} \Vvert_{\mathbf{w},p}^p
\leq \frac{1}{\mu_n \, c}\;\Phi_{\mu_n; g_n}(f^\star_{\mu_n;g_n})
\leq \frac{1}{\mu_n \, c}\; \Phi_{\mu_n;g_n}(\fs )\nonumber \\
&=&\frac{1}{\mu_n \, c}\left[\|Kf_o-g_n \|^2
+ \mu_n \Vvert \fs \Vvert^p_{\mathbf{w},p} \right]
\leq  \frac{1}{c} \left( \frac{\epsilon_n^2}{\mu_n}+\Vvert f^\dagger
\Vvert_{\mathbf{w},p}^p \right) ~,
\label{fmuubd}
\end{eqnarray}
where we have used, respectively,
the bound (\ref{bdL2Ban}), the fact that $f^\star_{\mu_n;g_n}$
minimizes $\Phi_{\mu_n;g_n}(f)$, $K\fs = Kf_o$
and the bound $\|Kf_o-g_n\|^2 \leq \epsilon_n^2$.
By the assumption \eref{mu-req}, $\epsilon_n^2 /\mu_n$ tends to zero for
$n\to\infty$ and hence can be bounded by a constant independent of $n$.
\newline
It
follows that the sequence $(f^\star_{\mu_n;g_n})_{_{n\in \N}}$ has at least
one weak
accumulation point, i.e. there exists a subsequence

$(f^\star_{\mu_{n_l};g_{n_l}})_{_{l \in {\N}}}$ that has a weak limit.
Because this sequence is bounded in the $\Vvert \; \Vvert $-norm,
by passing to a subsequence
 $\left( f^*_{\mu_{n_{l(k)}};g_{n_{l(k)}}}\right)_{k \in \N}$, we can
ensure that the $\Vvert f^*_{\mu_{n_{l(k)}};g_{n_{l(k)}}}
\Vvert_{\mathbf{w},p}$
constitute a converging sequence.
To simplify notation, we define $\widetilde\mu_k = \mu_{n_{l(k)}}$ and
${\widetilde f}_k =
f^\star_{\mu_{n_{l(k)}};g_{n_{l(k)}}}$; the $\widetilde{f}_k$ have the
same weak limit $\widetilde{f}$ as the $f^*_{\mu_{n_l};g_{n_l}}$.
We also define
$\widetilde{g}_k = g_{n_{l(k)}}$,
${\widetilde e}_k = \widetilde{g}_k-Kf_o$ and $\widetilde\epsilon_k =
\epsilon_{n_{l(k)}}$.
We shall show that $\widetilde{f}=\fs$.
\newline
Since each $\widetilde{f}_k $ is the minimizer of
$\Phi_{\widetilde\mu_k; \widetilde{g}_k}(f) $, by Proposition \ref{fix-min}, it
is a fixed point of the corresponding operator $\T$. Therefore, for any $\g \in
\Gamma$,
${\wf}_\g  =
\left<{\wf},\vpg \right>$ satisfies
\begin{eqnarray*}
{\wf}_\g & = &\lim_{k\to\infty}({\wf}_k)_\g  =
\lim_{k\to\infty} S_{\widetilde\mu_k \ag,p}[({\widetilde h}_k)_\g] \\
\hbox{\ \ with\ \ } ~~~  {\widetilde h}_k & = &
{\wf}_k + K^*(\widetilde g_k-K{\wf}_k)=
\wf_k + \K K (f_o - \wf_k)+ K^*{\widetilde e}_k ~.
\end{eqnarray*}
We now rewrite this as
\begin{equation}
\label{2terms}
{\wf}_\g=\lim_{k\to\infty}
\left(S_{\widetilde\mu_k \ag,p}[({\widetilde h}_k)_\g] -
({\widetilde h}_k)_\g \right) +
\lim_{k\to\infty} ({\widetilde h}_k)_\g ~ .
\end{equation}
By Lemma \ref{lm-4-1} the first limit in the right hand side is zero, since
$$
\left|S_{\widetilde\mu_k \ag,p}[({\widetilde h}_k)_\g] -
({\widetilde h}_k)_\g \right|
\leq  p\ \ag\,   \widetilde{\mu}_k
~| (\widetilde h_k)_\g |^{p-1} /2\leq p\ C\  \widetilde{\mu}_k
[ 3C + \widetilde\epsilon_k ] /2 \xrightarrow[k \to \infty]{~} 0~,
$$
where we have used $\|K\|<1$ ($C$ is some constant depending on $\ag$). Because
$\lim_{k
\to
\infty}\Vert\widetilde{e}_k\Vert=0$, and {\em w}$\,$--$\lim_{k \to \infty}
\wf_k =
\wf$,  it then follows from \eref{2terms} that
\begin{equation*}
{\wf}_\g = \lim_{k\to\infty} ({\widetilde h}_k)_\g = {\wf}_\g +
[\K K(\fs -{\wf})]_\g \ .
\end{equation*}
Since this holds for all $\g$, it follows that
$\K K(\fs-{\wf})=0$. If $\mN(K)=\{0\}$, then this allows us
immediately to conclude that $\wf=\fs$. When $\mN(K) \neq \{0\}$,
we can only conclude that $\fs -\wf \in \mN(K)$. Because $\fs$
has the smallest $\Vvert ~ \Vvert _{\mathbf{w},p}$--norm among all
$f \in \mathcal{S}=\{f; Kf=Kf_o\}$, it follows that
$\Vvert \wf \Vvert _{\mathbf{w},p} \geq \Vvert \fs \Vvert _{\mathbf{w},p}$.
On the other hand, because
the ${\wf}_k$ weakly converge to ${\wf}$, and therefore, for all $\g$,
$({\wf}_k)_\g \to {\wf}_\g $ as $k\to\infty$,  we
can use Fatou's lemma to obtain
\begin{equation}
\Vvert {\wf} \Vvert^p_{\mathbf{w},p} =
{\sum_\g} \ag|{\wf}_\g |^p  \leq
\limsup_{k \to\infty} {\sum_\g} \ag |({\wf}_k)_\g|^p
= \limsup_{k \to\infty}\Vvert {\wf}_k\Vvert_{\mathbf{w},p}^p
= \lim_{k \to\infty}\Vvert {\wf}_k\Vvert_{\mathbf{w},p}^p \ .
\label{now123}
\end{equation}
It then follows from (\ref{fmuubd}) that
\begin{equation}
\lim_{k \to\infty}\Vvert {\wf}_k\Vvert_{\mathbf{w},p}^p
\leq \lim_{k \to\infty}\left[\frac{\widetilde\epsilon_k^2}{ \widetilde\mu_k}
+\Vvert \fs \Vvert_{\mathbf{w},p}^p\right] = \Vvert \fs
\Vvert_{\mathbf{w},p}^p
\leq \Vvert {\wf}\Vvert_{\mathbf{w},p}^p \ .
\label{now124}
\end{equation}
Together, the inequalities (\ref{now123}) and (\ref{now124}) imply that
\begin{equation}
\lim_{k \to\infty}\Vvert {\wf}_k\Vvert_{\mathbf{w},p}
= \Vvert \fs \Vvert_{\mathbf{w},p}= \Vvert {\wf}\Vvert_{\mathbf{w},p} \ .
\label{now125}
\end{equation}
Since $\fs$ is the unique element in $\mathcal{S}$ of minimal
$\Vvert ~ \Vvert_{\mathbf{w},p}$--norm, it follows
that
${\wf}=\fs$.
The same argument holds for any other weakly converging subsequence of
$(f^\star_{\mu_n;g_n})_{n \in {\N}}$; it follows that the sequence
itself converges weakly to $\fs$.
Similarly we conclude from \eref{now125} that $\lim_{n \to \infty}
\Vvert f^\star_{\mu_n;g_n} \Vvert_{\mathbf{w},p} =
\Vvert \fs _{\mu_n;g_n} \Vvert_{\mathbf{w},p}~$.
It then follows from Lemma \ref{lm-conv} that the $f^\star_{\mu_n;g_n}$
converge to $\fs$ in the $\cH$-norm. \hfill \QED

\bigskip

\begin{remark} {\rm Even when $p=1$ and $N(K) \neq \{0\}$, it may still be the
case that, for any $f_o \in \cH$, there is a unique element $\fs$ of minimal
norm in $\mathcal{S}=\{f \in \cH; Kf=Kf_o \}$. (For instance, if
$K$ is diagonal in the $\varphi_\g$--basis, with some zero eigenvalues,
then the unique minimizer $\fs$ in $\mathcal{S}$ is given by setting to zero
all the components of $f_o$ corresponding to $\g$ for which $K \vpg = 0$.)
In this case the proof still applies, and we still have norm--convergence
of the $f^\star_{\mu(\epsilon),\mathbf{w},p;g}$ to $\fs$ if $\mu(\epsilon)$
satisfies \eref{mu-req} and $\|g -Kf_o\| \leq \epsilon \to 0$.}
\end{remark}

\subsection{Stability estimates}

The regularization theorem of the previous subsection
gives no information on the rate at
which the regularized solution approaches the exact solution when the noise
(as measured by $\epsilon$) decreases to
zero. Such rates are not available in the general case, but can be derived
under additional assumptions, discussed below.
For the remainder of this section we shall assume
that the operator $K$ is invertible on its range, i.e. that
$\mN(K)=\{0\}$. Suppose that the
unknown exact solution of the problem, $f_o$, satisfies the constraint
$\Vvert f_o \Vvert_{\mathbf{w},p} \leq \rho$, where $\rho>0$ is given; in
other words, we know a priori
that the unknown solution lies in the ball around the origin with radius $\rho$
in the Banach space $\cB_{\mathbf{w},p}$; we shall denote this ball
by $\rm{B}_{\mathbf{w},p}(0,\rho)$. If we also know that $g$ lies within
a distance $\epsilon$ of $Kf_o$ in $\cH'$, then
we can localize the exact
solution within the set
\begin{equation*}
{\cal F}(\epsilon,\rho) =
\{ f\in \cH ; \, \| Kf-g \| \leq \epsilon\, ,
\, \Vvert f\Vvert_{\mathbf{w},p} \leq \rho \}\ .
\end{equation*}
The diameter of this set is a measure of the uncertainty of the
solution for a given prior and a given noise level $\epsilon$. The maximum
diameter of ${\cal F}$, namely diam(${\cal F}$)=$\sup\{ \| f-f'\|;\
f,f' \in {\cal F}\}$ is bounded by $2M(\epsilon,\rho)$, where
$M(\epsilon,\rho)$, defined by
\begin{equation}
M(\epsilon,\rho)=\sup\{\| h\|;\, \| Kh \|
\leq \epsilon \, ,\, \Vvert h\Vvert_{\mathbf{w},p} \leq \rho\} \ ,
\label{MC}
\end{equation}
is called the {\it modulus of continuity} of $K^{-1}$ under the
prior. (We have once more dropped the explicit reference in our
notation to the dependence on $\mathbf{w}$ and $p$.)
If \eref{compact-emb} is satisfied, then
the ball $\rm{B}_{\mathbf{w},p}(0,\rho)$ is compact in $\cH$, and it
follows from a general topological lemma
(see e.g. \cite{Eng96}) that $M(\epsilon,\rho) \to 0$
when $\epsilon \to 0~$; the uncertainty on the solution
thus vanishes in this limit. However,
this topological argument, which holds for any regularization
method enforcing the prior $\Vvert f_o \Vvert_{\mathbf{w},p} \leq \rho$,
does not tell us anything about the rate of convergence
of any specific method.

In what follows, we shall systematically assume that \eref{compact-emb}
is satisfied. We shall also make additional assumptions that will make it
possible to derive more precise convergence results.
Our specific regularization method consists in taking the
minimizer $f^*_{\mu;g}$ of the functional
$\Phi_{\mu; g}(f )$ given by
(\ref{phimu}) as
an estimate of the exact solution $f_o$, where we
leave any links between $\mu$ and $\epsilon$ unspecified for
the moment. (Because of the compactness argument above, we could
conceivably dispense with \eref{mu-req}; see below.) An upper
bound on the reconstruction error $\| f^*_{\mu;g} - f_o \|$ ,
valid for all $g$ such that $\|g-Kf_o\| \le \epsilon$, as well as uniformly
in $f_o$, is
given by the following {\it modulus of convergence}:
\begin{equation}
M_\mu(\epsilon,\rho) = \sup\{ \| f^*_{\mu;g} -f\|;\ f \in \cH,\, g \in \cH'
\, ,
\| Kf - g \| \leq \epsilon\, ,\, \Vvert f \Vvert_{\mathbf{w},p} \leq \rho \}\ .
\label{modcv}
\end{equation}
The decay of this modulus of convergence as $\epsilon \to 0$ is governed by the
decay of the modulus of continuity \eref{MC}, as shown by the following
proposition.
\begin{proposition}
\label{bestcv}
The modulus of convergence {\rm \eref{modcv}} satisfies
\begin{equation}
M(\epsilon,\rho) \leq M_\mu(\epsilon,\rho) \leq M(\epsilon
+\epsilon',\rho + \rho')\ .
\label{stability}
\end{equation}
where
\begin{equation}
\epsilon' = \left(\epsilon^2 +\mu \rho^p\right)^{\frac{1}{2}}
~~~~\mbox{and}~~~~
\rho' = \left(\rho^p + \frac{\epsilon^2}{\mu}\right)^{\frac{1}{p}}\ .
\label{primes}
\end{equation}
and $M(\epsilon,\rho)$ is defined by {\rm \eref{MC}}\ .
\end{proposition}
{\em Proof:} We first note
that $\Phi_{\mu ; g}(f^*_{\mu;g}) \leq\Phi_{\mu;g}(f_o) \leq \epsilon^2 + \mu
\rho^p$ because $f^*_{\mu;g}$ is the minimizer of $\Phi_{\mu ; g}(f)$
and $f_o \in {\cal F}(\epsilon, \rho)$.
It follows that
$$
\| Kf^*_{\mu;g} - g\|^2 \leq \Phi_{\mu; g}(f^*_{\mu;g}) \le
\epsilon^2 + \mu \rho^p  ~~\mbox{and}~~
\mu  \Vvert f^*_{\mu;g}\Vvert_{\mathbf{w},p}^p \leq \Phi_{\mu; g}(f^*_{\mu;g} )
 \le \epsilon^2 + \mu \rho^p
$$
or, equivalently, $f^*_{\mu;g} \in {\cal F}(\epsilon',\rho')$ with
$\epsilon'$ and $\rho'$ given by \eref{primes}.
The modulus of convergence (\ref{modcv}) can then be bounded as follows,
using the
triangle inequality. Indeed, for any $f \in {\cal F}(\epsilon,\rho)$ and
$f' \in {\cal F}(\epsilon',\rho')$, we have
$
\| K(f - f') \|
\leq \epsilon+\epsilon'
$
and
$
\Vvert f - f' \Vvert_{\mathbf{w},p} \leq \rho+\rho'\ .
$
and we immediately obtain from the definition of (\ref{MC}) the
upper bound in (\ref{stability}). To derive the lower bound, observe that
for the
particular choice $g=0$,
the minimizer $f^*_{\mu;g}$ of the functional (\ref{phimu}) is
$f^*_{\mu;0}=0$. The desired lower bound then follows immediately upon
inspection
of the two definitions \eref{MC} and \eref{modcv}. \hfill\QED

\bigskip

Let us briefly discuss the meaning of the previous proposition.
The modulus of continuity $M(\epsilon,\rho)$ yields the best possible
convergence rate for any
regularization method that enforces the error bound and the prior
constraint defined by \eref{MC}.
Proposition \ref{bestcv} provides a relation between the modulus of
continuity and  the convergence rate $M_\mu(\epsilon,\rho)$ of the
specific regularization method considered in this paper, which
is defined by the minimization of the functional \eref{phimu}. Optimizing
the upper bound
in  \eref{stability} suggests the choice
$\mu=\epsilon^2/\rho^p$, yielding
$\epsilon'=\sqrt{2}\;\epsilon\ $ and $\rho'=2^{1/p}\rho\ $.  With these
choices,
we ensure that $f^*_{\mu;g} \rightarrow f_o$ when $\epsilon
\rightarrow 0$, i.e.
that the problem is {\it regularized}, provided we can show
that the modulus of continuity
tends to zero with $\epsilon$.
 Moreover, once we establish its rate of
decay (see below), we know that our regularization method is (nearly)
optimal in the sense that the modulus of convergence (\ref{modcv}) will decay
{\em at the same rate as the optimal rate} given by the modulus of stability
$M(\epsilon,\rho)$\ (We call it {\em nearly} optimal because, although the
rate of
decay is optimal, the constant multiplier probably is not.)
Note that because of the assumption of compactness of
the ball $\rm{B}_{\mathbf{w},p}(0,\rho)$ (which amounts to assuming that
\eref{compact-emb} is satisfied), we achieve regularization
even in some cases where $\epsilon^2 / \mu$
does not tend to zero for $\epsilon \to 0$, which is a case not covered by
Theorem \ref{regthm}.

In order to derive upper or lower bounds on $M(\epsilon,\rho)$, we must
know more information about the operator $K$. The following proposition
illustrates how such information can be used.

\begin{proposition}
\label{convrate-gen}
Suppose that there exist sequences $\mathbf{b}=(b_\g)_{\g \in \Gamma}$
and $\mathbf{B}=(B_\g)_{\g \in \Gamma}$ satisfying
$\forall \,\g \in \Gamma ~:~ 0 < b_\g,
\, B_\g < \infty\,$ and such that, for all $h$ in $\cH$,
\begin{equation}
\label{extraK}
\sum_{\g \in \Gamma} b_\g |h_\g|^2 \leq \|K h \|^2
\leq \sum_{\g \in \Gamma} B_\g |h_\g|^2 ~.
\end{equation}
Then the following upper and lower bounds hold for $M(\epsilon,\rho)$:
\begin{eqnarray}
\label{lowerbound}
M(\epsilon,\rho) &\geq & \max_{\g \in \Gamma}
\left[ \min \left(\rho \ag^{-1/p},\epsilon B_\g^{-1/2} \right) \right]~, \\
\label{upperbound}
M(\epsilon,\rho) &\leq & \min_{\Gamma = \Gamma_{_{\!1}} \cup \Gamma_{_{\!2}}}
\sqrt{ \frac{\epsilon^2}{\min_{\g \in \Gamma_{_{\!1}}}b_\g} +
\frac{\rho^2}{\min_{\g \in \Gamma_{_{\!2}}}\ag^{2/p}} } ~.
\end{eqnarray}
\end{proposition}
{\em Proof:}
To prove the lower bound, we need only exhibit one particular $h$ such that
$\|Kh\| \le \epsilon$ and $\Vvert h \Vvert_{\mathbf{w},p}\le \rho$, for which
$\|h\|$ is given by the right hand side of
\eref{lowerbound}. For this we need only identify the index ${\g_{_m}}$ such
that, $\forall \g \in \Gamma \,$,
$$
\nu = \min \left(\rho w_{\g_{_m}}^{-1/p},\epsilon B_{\g_{_m}}^{-1/2} \right)
\ge  \min \left(\rho \ag^{-1/p},\epsilon B_\g^{-1/2} \right) ~,
$$
and choose $h = \nu \varphi_{\g_{_m}}$. Then $\Vvert h
\Vvert_{\mathbf{w},p} =
\nu \, w_{\g_{_m}}^{1/p} \le \rho$ and $\|K h \| \le  \nu B_{\g_{_m}}^{1/2}
\le \epsilon ~$; on the other hand, $\nu$ equals the right hand side of
\eref{lowerbound}.

On the other hand, for any partition of $\Gamma$ into two subsets,
$\Gamma = \Gamma_{_{\!1}} \cup \Gamma_{_{\!2}}$,
and any $h \in \{h; \|Ku\|\le \epsilon, ~\Vvert u \Vvert_{\mathbf{w},p} \le
\rho \}$, we have
\begin{eqnarray*}
\sum_{\g \in \Gamma}|h_\g|^2 &=& \sum_{\g \in \Gamma_{_{\!1}}} |h_\g|^2 +
\sum_{\g \in \Gamma_{_{\!2}}} |h_\g|^2 \\
&\le & \max_{\g' \in \Gamma_{_{\!1}}} [b_{\g'}^{-1} ]
\sum_{\g \in \Gamma_{_{\!1}}} b_\g |h_\g|^2 +
\max_{\g' \in \Gamma_{_{\!2}}}[ w_{\g'}^{-2/p}][\max_{\g'' \in
\Gamma_{_{\!2}}} w_{\g''}
|h_{\g''}|^p]^{\frac{2}{p}-1 } \sum_{\g \in \Gamma_{_{\!2}}} \ag |h_\g|^p \\
& \le & \max_{\g' \in \Gamma_{_{\!1}}} [b_{\g'}^{-1}] \epsilon^2
+ \max_{\g' \in \Gamma_{_{\!2}}}[ w_{\g'}^{-2/p}] \rho^2 ~.
\end{eqnarray*}
Since this is true for any partition $\Gamma = \Gamma_{_{\!1}} \cup
\Gamma_{_{\!2}}$, we
still have an upper bound, uniformly valid for all $h \in
\{u; \|Ku\|\le \epsilon, ~\Vvert u \Vvert_{\mathbf{w},p} \le \rho \}$, if we
take the minimum over all such partitions. The upper bound on
$M(\epsilon,\rho)$ then follows upon taking the square root.
\hfill \QED

\bigskip

To illustrate how Proposition \ref{convrate-gen} could be used, let us apply
it to one particular example, in which we choose
the $(\vpg)$--basis with respect to which
the $\Vvert~\Vvert_{\mathbf{w},p}$--norm is defined to be a wavelet basis
$(\Psi_\lambda)_{\lambda \in \Lambda}$. As
already pointed out in subsection 1.4.1
the Besov spaces $B^s_{p,p}(\R^d)$ can then be identified with the Banach
spaces
$\mathcal{B}_{\mathbf{w},p}$ for the particular choice $w_\lambda =
2^{\sigma  p
|\lambda|}$, where $\sigma = s + d\left( \frac{1}{2}-\frac{1}{p} \right)$
is assumed to be non-negative. For $f \in B^s_{p,p}(\R^d)$, the Banach norm
$\Vvert f \Vvert_{\mathbf{w},p}$ then coincides with the Besov norm
$\VVert f \VVert_{s,p} = \left[ \sum_{\lambda \in \Lambda} w_{\lambda}
|\left<f,\Psi_\lambda \right>|^p \right]^{1/p} ~.$
Let us now consider an inverse problem for the operator $K$ with such a
Besov prior.
If we assume that the operator
$K$ has particular smoothing properties, then we can
use these to derive bounds on the corresponding modulus of continuity, and thus
also on the rate of convergence for our regularization algorithm.
In particular, let us assume that the operator $K$ is a smoothing
operator of order
$\alpha$, a property which can be formulated as an equivalence between the norm
$\| Kh\|$ and the norm of $h$ in a Sobolev space of
negative order $H^{-\alpha}$, i.e. in a Besov space $B^{-\alpha}_{2,2}$
(see e.g.
\cite{Eng96}, \cite{Lou97}, \cite{DeV98}, \cite{Coh02}). In other words,
we assume that
for some $\alpha >0$, there exist constants  $A_{\ell}$ and $A_u$,  such that,
for all $h \in L^2(\R^d)$,
\begin{equation}
\label{ell}
A_{\ell}^2 \sum_\lambda 2^{-2\vert\lambda\vert\alpha}\; \vert
h_\lambda\vert^2 \leq \| K h\|^2 \leq A_u^2\sum_\lambda
2^{-2\vert\lambda\vert\alpha}\;
\vert h_\lambda\vert^2\ .
\end{equation}
The decay rate of the modulus of continuity is then characterized as follows.
\begin{proposition}
If the operator $K$ satisfies the smoothness condition {\rm(\ref{ell})},
then the
modulus of continuity $M(\epsilon,\rho)$, defined by
$$
M(\epsilon,\rho)= \max\{ \|h\|; \|K h\| \leq \epsilon,
\VVert h \VVert_{s,p} \le \rho \} ~,
$$
satisfies
\begin{equation}
\label{M}
c \left( \frac{\epsilon}{A_u} \right)^{\frac{\sigma}{\sigma+\alpha}}
\rho^{\frac{\alpha}{\sigma+\alpha}} \leq M(\epsilon,\rho) \leq C
\left( \frac{\epsilon}{A_\ell}\right)^{\frac{\sigma}{\sigma+\alpha}}
\rho^{\frac{\alpha}{\sigma+\alpha}}
\end{equation}
where $\sigma = s + d\left( \frac{1}{2}-\frac{1}{p} \right)\geq0$,
and $c$ and $C$ are constants depending on $\sigma$
and $\alpha$ only.
\end{proposition}
{\em Proof:}
By \eref{ell}, the operator $K$ satisfies
\eref{extraK} with $b_\lambda=A_{\ell}^2\, 2^{-2\vert\lambda\vert\alpha}$
and $B_\lambda=A_u^2\, 2^{-2\vert\lambda\vert\alpha} \,$.
\newline
It then follows from \eref{lowerbound} that
$$
M(\epsilon, \rho) \geq \max_\lambda \left[ \min \left(
\rho \, 2^{- \sigma |\lambda|} , \frac{\epsilon}{A_u}\, 2^{\alpha |\lambda|}
\right) \right] ~;
$$
if $x=|\lambda|$ could take on all positive real values, then one easily
computes
that this max-min would be given for
$x= - [\log_2(\epsilon/\rho A_u)]/(\alpha + \sigma)$, and would be equal to
$(\epsilon/A_u)^{\sigma/(\alpha + \sigma)} \rho^{\alpha/(\alpha + \sigma)}$.
Because $|\lambda|$ is constrained to take only the values in $\N$, the max-min
is guaranteed only to be within a constant of this bound (corresponding
to an integer neighbor of the optimal $x$), which leads
to the lower bound in \eref{M}.

For the upper bound \eref{upperbound}, we must partition the index set.
Splitting $\Lambda$ into $\Lambda_{_1}=\{\lambda; |\lambda| < J\}$ and
$\Lambda_{_2}=\{\lambda; |\lambda| \ge J\}$, we find that
$$
\frac{\epsilon^2}{\min_{\lambda \in \Lambda_{_{1}}}b_\lambda} +
\frac{\rho^2}{\min_{\lambda \in \Lambda_{_{2}}}w_\lambda^{2/p}}
=\frac{\epsilon^2}{A_\ell^2}\, 2^{2\alpha(J-1)} +  \rho^2\, 2^{-2 \sigma J} ~.
$$
The minimizing partition for $\Lambda$ thus corresponds with the minimizing
$J$ for the right hand side of this expression. This value
for  $J$ is an integer neighbor of
$y= - [\log_2(\epsilon/\rho A_\ell)]/(\alpha + \sigma)$, which leads to the
upper bound in \eref{M}.
\hfill \QED

\bigskip

The stability estimates we have derived are standard in regularization
theory for the special case $p=2$; they were first extended
to the case $1\leq p < 2$
in \cite{DeM02}. They show the interplay between the smoothing order of the
operator
characterized by $\alpha$ and the assumed smoothness of
the solutions characterized by $\sigma=s+d(\frac{1}{2}-\frac{1}{p})$
(for Besov spaces, we recall that this amounts to
solutions
having $s$ derivatives in $L^p$). For
$\sigma/(\sigma+\alpha)$ close to one, the problem is mildly ill-posed,
whereas the stability degrades for large $\alpha$. Note that if
the bound (\ref{ell}) were replaced by another one,
in which the decay of the $b_\lambda$ and
$B_\lambda$ was given by an exponential decay in $D=2^{|\lambda|}$ (instead
of the much slower decaying negative power $D^{-2\alpha}$) of \eref{ell},
then the modulus of continuity would tend to zero only as an
inverse power
of $| \log\epsilon|$. This is the so-called {\it logarithmic continuity}
which has been extensively discussed in the case $p=2$, and which extends,
as shown
by an easy application of Proposition \ref{convrate-gen}, to $1 \le p <2$.

\section{An illustration}
We provide a simple illustration of the behavior of the algorithms based on
minimizing
$\Phi_{\mu \mathbf{w}_0,1}$ and
$\Phi_{\mu \mathbf{w}_0,2}$ for a two-dimensional deconvolution
problem,
considering a class of objects consisting of small bright sources on a dark
background. The image is discretized on a $256 \times 256$ array,
denoted by  $f$. The convolution operator $K$ is
implemented by multiplying the discrete Fourier transform (DFT) of $f$
by a low-pass, radially symmetrical filter and then
taking the inverse DFT to obtain the data $g=Kf$ (data were zero padded on
a larger
$512 \times 512$ array when taking DFTs). The filter was equal
(in the Fourier domain) to the convolution with itself
of the characteristic function of a disk with radius equal to $0.1$ times the
maximum frequency determined by the image grid sampling;
this filter provides a discrete model of a diffraction-limited imaging
system with
incoherent light. Pseudo-random Poisson noise was
added to the data array $g$, for a total number of $10000$ photons,
corresponding to about $25$ photons for the data pixel with the maximum
intensity.

\begin{figure}[h!tb]
\begin{center}
\epsfig{figure=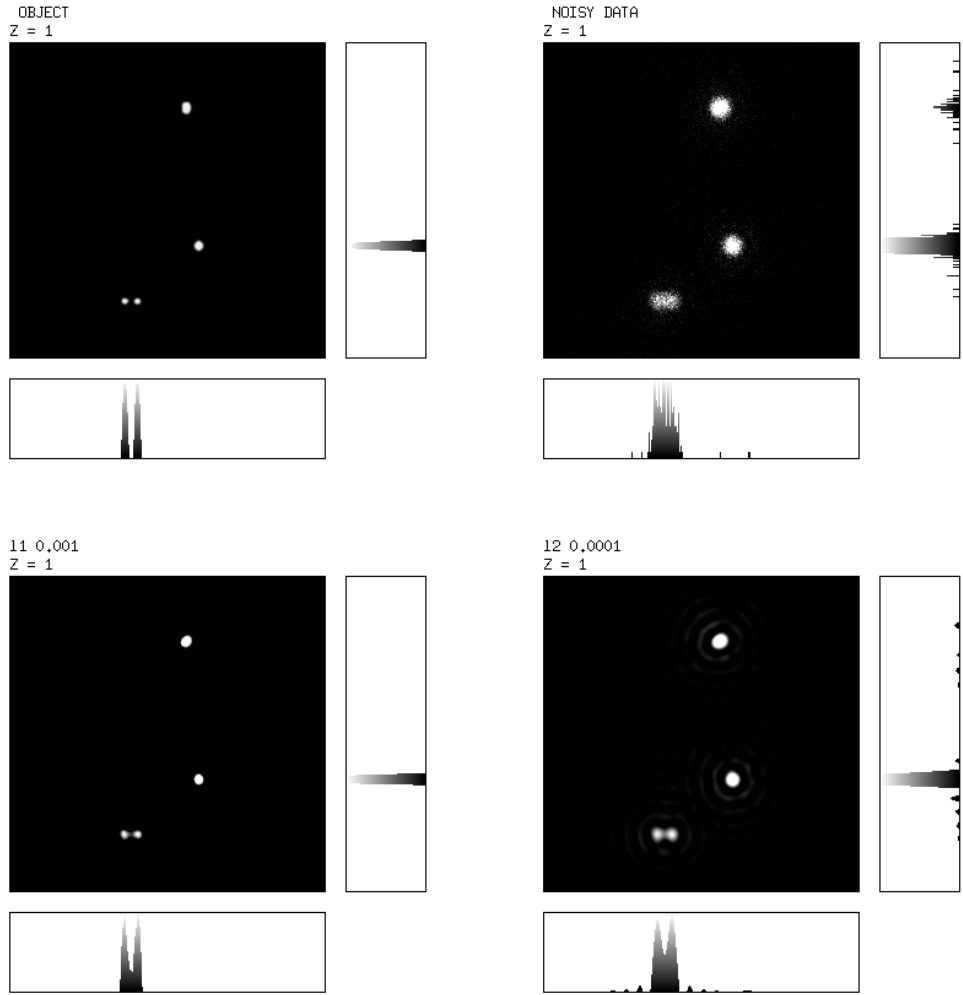, width= 5 in}
\caption{The object $f$ (top left), the image $g$ after convolving
with a radially symmetrical low-pass filter and adding pseudo-random
Poisson noise (top right), and the minimizers of $\Phi_{\mu_1 \mathbf{w}_0,1}$
(bottom left) and $\Phi_{\mu_2 \mathbf{w}_0,2}$ (bottom right). The values
$\mu_1=0.001$ and $\mu_2=0.0001$ have been selected separately for the
$\ell^1$- and $\ell^2$-cases, to obtain a balance between
sharpness and ringing and noise. }
\end{center}
\end{figure}

The top of Figure 1 shows the object $f$ (four ellipses of axis
$7.5$ or $5.0$ pixels, slightly smoothed to avoid blocking effects)
and the data $g=Kf$. The figure also shows intensity
distributions along two lines in the object and data arrays; along the
horizontal line we see
how two close sources in $f$ give rise to a joint blur in $g$.
The bottom of Figure 1 shows the reconstructions obtained after
2000 steps of the iterative thresholded Landweber algorithm, which
accurately approximate the minimizers of $\Phi_{\mu_1 \mathbf{w}_0,1}$ and
$\Phi_{\mu_2 \mathbf{w}_0,2}$. The parameters $\mu_1$ and $\mu_2$
are selected separately for each case, in order to achieve a good balance
between
sharpness and ringing and noise.

\begin{figure}[h!tb]
\begin{center}
\epsfig{figure=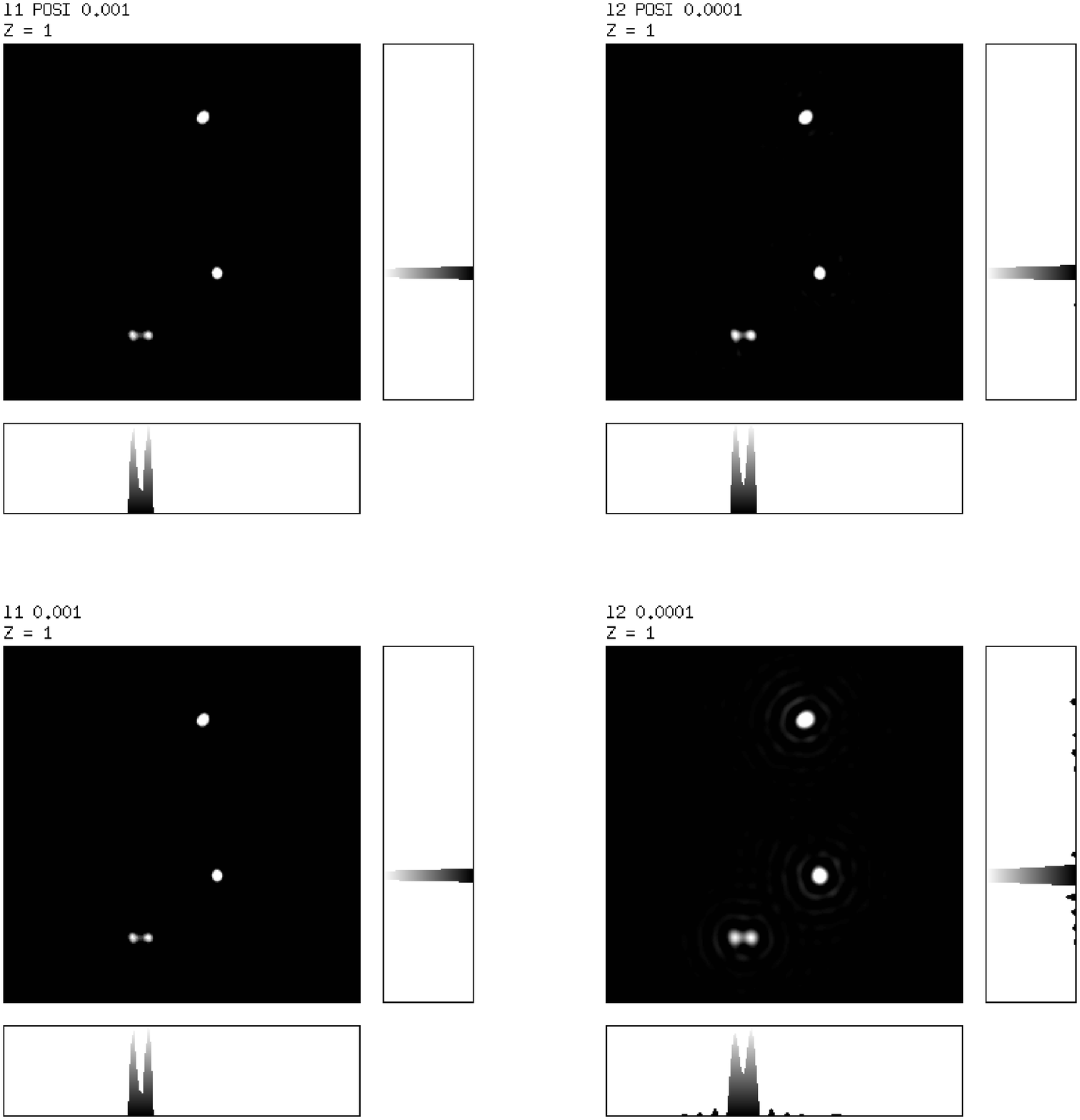, width= 5 in}
\caption{A comparison to illustrate the impact of the positivity constraint,
imposed at every iteration step. On the left are the fixed points for
$P_{\cal C} \mathbf{S}_{\mu_1 \mathbf{w}_0,1}$ (top) and
$\mathbf{S}_{\mu_1 \mathbf{w}_0,1}$ (bottom); on the right
those of $P_{\cal C} \mathbf{S}_{\mu_2 \mathbf{w}_0,2}$ (top) and
$\mathbf{S}_{\mu_2 \mathbf{w}_0,2}$ (bottom). The data
and the values of $\mu_1,~\mu_2$ are the same as in Figure 1.}
\end{center}
\end{figure}

As expected for an example of this type, the minimizer of
$\Phi_{\mu_1 \mathbf{w}_0,1}$ does a better job at
resolving the two close sources on the horizontal line; it also gives a better
concentration of the source lying along the vertical line.
Because the object $f$ is positive, we can
apply Remark 3.12  and use $P_{\cal C} S_{\mu \mathbf{w}_0,p}$ instead
of $S_{\mu \mathbf{w}_0,p}$, where $P_{\cal C}$ is the projection
onto the  convex cone of $256 \times 256$ arrays that take only non-negative
values.

The results are
shown in Figure 2; on the top is the $2000$-th iterate for the case with
$P_{\cal C}$,  with the case $p=1$
(left) and $p=2$ (right). In each case we used the same values of
$\mu_1$ and $\mu_2$ as for
the reconstructions without positivity constraint, which are shown for
comparison on the bottom of
Figure 2. Exploiting the positivity constraint leads to better resolution and
less ringing for this example, where the background is zero.

\begin{figure}[h!tb]
\begin{center}
\epsfig{figure=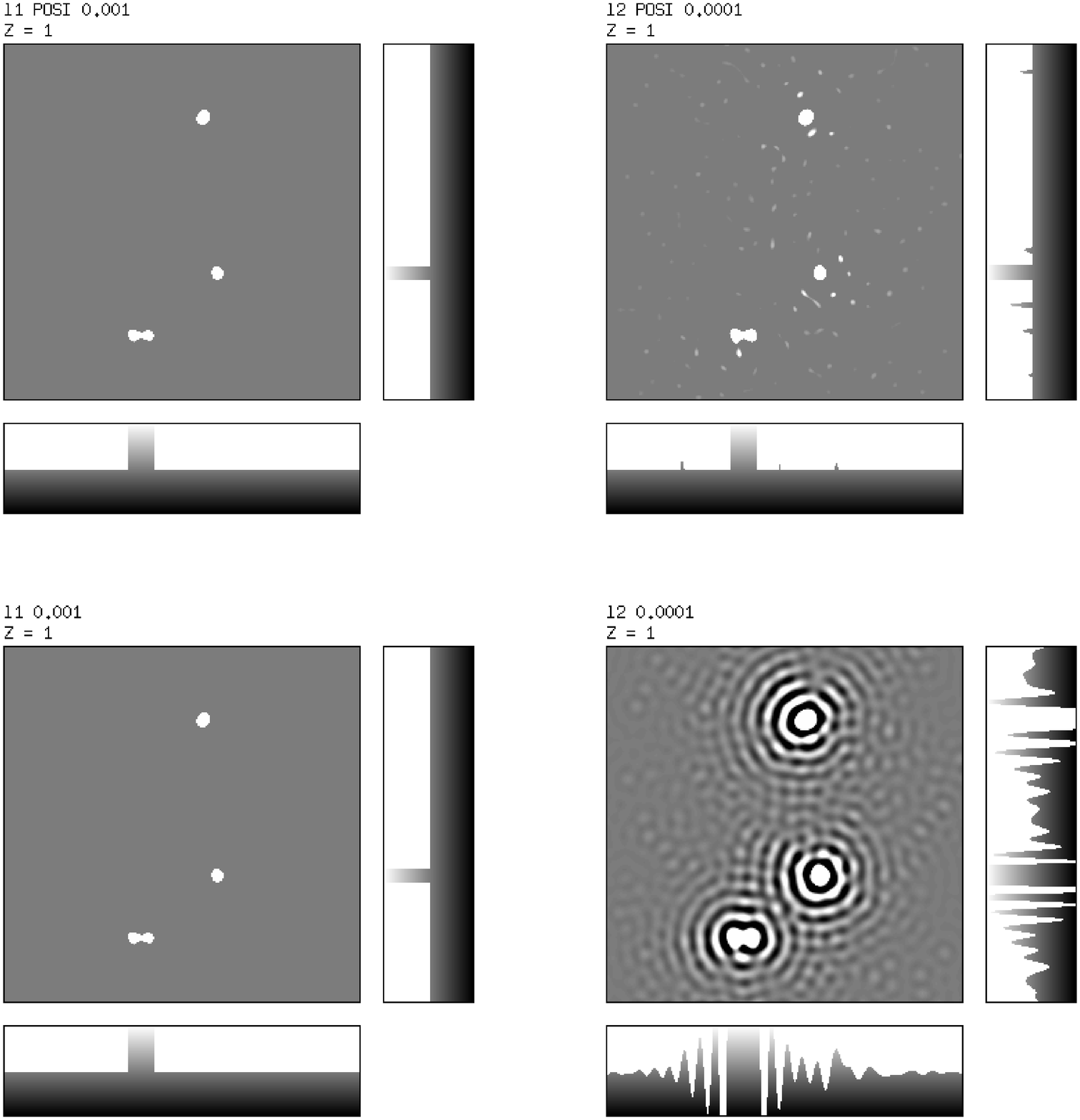, width= 5 in}
\caption
{A different view of the four solutions in Figure 2, with a different
dynamic range for the image intensity gray scale, to highlight ringing
and other artifacts.}
\end{center}
\end{figure}

Figure 3 gives a different view of the same solutions, with
a compressed gray scale ranging from $- 2 \%$ (darker) to $+ 2 \%$ (lighter)
of the maximum intensity in the original object. This has the effect of
highlighting
the ringing effects and the noise. Both
ringing and noise are seen to be less pronounced for the minimizer of
$\Phi_{\mu_1 \mathbf{w}_0,1}$ (top left)
than that of $\Phi_{\mu_2 \mathbf{w}_0,2}$. Although the introduction of
the positivity constraint removes the
ringing phenomenon (top of Figure 3), we nevertheless see that noise is better
suppressed with $p=1$.

To produce Figures 1 to 3, the same program was used in every case; the only
change was the choice
of the nonlinear operator applied at the $n$-th iteration step to
$f^{n-1} + K^* (g - K f^{n-1})$.
For realistic applications on data of this type, more sophisticated algorithms
exist. With the
$\ell^2$--penalty, for instance, the reconstructions in our simple example
can be
obtained directly by a regularized Fourier deconvolution.
These examples are included to illustrate the differences
that can be achieved by the choice of $p$, and do not constitute a claim
that the
iterative algorithms discussed
in this paper are optimal. The ``data'' in this example are also only
simple-minded
caricatures of quasi
point-sources data sets. While similar examples may have applications in
astronomy,
most natural
images have a much richer structure. However, as is abundantly documented, the
wavelet transforms of
natural images tend to have distributions that are sparse. A similar
improvement in
accuracy can be expected by applying $\ell^1$ rather than $\ell^2$
penalizations on the
wavelet coefficients
in inverse problems involving natural images, similar to the gain achieved in
denoising with
a soft thresholding rather than with a quadratic penalty.

\section{Generalizations and additional comments}

The algorithm proposed in this paper can be generalized in several directions,
some of which we list here, with brief comments.

The penalization functionals $\Vvert f \Vvert_{\mathbf{w},p}$ we have used are
symmetric, i.e. they are invariant under the exchange of $f$ for $-f$. We can
equally well consider penalization functionals that treat positive and
negative
values of the $\fg$ differently. If $(w^+_\gamma)_{\gamma \in \Gamma}$ and
$(w^-_\gamma)_{\gamma \in \Gamma}$ are two sequences of
strictly positive numbers,
then we can consider the problem of minimizing the functional
\begin{equation}
\label{asfunct}
\Phi_{\mathbf{w}^+, \mathbf{w}^-,p} (f)= \Vert Kf-g\Vert^2 +
\sum_{\gamma \in \Gamma}
((w^+_\gamma) [\fg]_+^p + (w^-_\gamma) [\fg]_-^p)
\end{equation}
where, for $r \in \mathbb{R},~ r_+ = \max(0,r),~ r_- = \max(0,-r)$. One easily
checks that all the arguments in this paper can be applied equally well (after
some straightforward modifications) to the general functional \eref{asfunct},
provided we replace the thresholding functions $S_{w_\gamma,p}$ in the
iterative
algorithm by $S_{w^+_\gamma, w^-_\gamma, p}$, where, for $p>1$,
\begin{equation*}
S_{w^+, w^-, p} = \left( F_{w^+, w^-, p}\right)^{-1}\hbox{\ \ with\ \ }
F_{w^+, w^-, p}(x) = x +\frac{p}{2}\, w^+ [x]_+^{p-1} -
\frac{p}{2}\, w^- [x]_-^{p-1} ~,
\end{equation*}
and for $p=1$,
\begin{equation*}
S_{w^+, w^-, 1} = \left\{
 \begin{array}{ccl} x + w^-/2 ~&~ \mbox{if} ~& x \leq - w^-/2 \\
0 ~&~ \mbox{if} ~& - w^-/2 < x < w^+/2
\\ x- w^+/2 ~&~ \mbox{if} ~& x \geq  w^+/2.
\end{array} \right.
\end{equation*}

The above applies when the $\fg$ are all real; a generalization to complex
$\fg$
is not straightforward. When dealing with complex functions, one could
generalize
the penalization $\sum_{\gamma \in \Gamma} w_\gamma |\fg|^p$ to
$\sum_{\gamma \in \Gamma, |\fg| \neq 0} w_\gamma({\rm arg} \fg) |\fg|^p$,
where the weight coefficients have been replaced by strictly positive
$2\pi$--periodic $C^1$--functions on the $1$--torus $\mathbb{T} =
\{x \in \mathbb{C}, |x|=1\}$.
It turns out, however, that the variational
equation for $e^{i \arg\fg}=\fg |\fg|^{-1}$ then no longer uncouples
from that for $|\fg|$ (as it does in the case where $w_\gamma$ is a constant),
leading to a more complicated ``generalized thresholding'' operation in
which the
absolute value and phase of the complex number $S_{w,p}(\fg)$ are given
by a system of
coupled nonlinear equations.

When the $(\varphi_\gamma)_{\gamma \in \Gamma}$--basis is chosen to be a
wavelet basis,
then we saw in subsection 1.4.1 that is is possible to make the
$\Vvert ~ \Vvert_{\mathbf{w},p}$--norm equivalent to the Besov-norm
$\VVert ~ \VVert_{s,p}$, by choosing the weight for
$| \langle f, \Psi_\lambda \rangle |^p$ to be given by
$w_\lambda = 2^{|\lambda| \sigma p}$, where $|\lambda |$ is the scale of
wavelet
$\Psi_\lambda$. The label $\lambda$ contains much more information than just
the scale, however, since it also indicates the location of the wavelet,
as well as
its ``species'' (i.e. exactly which combination of $1$-dimensional scaling
functions
and wavelets is used to construct the product function $\Psi_\lambda$).
One could
choose the $w_\lambda$ so that certain regions in space are given
extra weight, or
on the contrary de-emphasized, depending on prior information.
In pixel space, prior information on the support of the object to be
reconstructed can be easily enforced by simply
setting the
corresponding weights to very small values,
or by choosing very large weights outside
the object support. This type of constraint is of uttermost importance
to achieve superresolution in inverse problems in optics and imaging
(see e.g. \cite{Ber96}).
 When thresholding in the wavelet domain,
a constraint on the object support can be enforced in a similar way due to the
good spatial localization of the wavelets.
If no a priori information is known,
one could
even imagine repeating the wavelet thresholding
algorithm several times, adapting
the weights $w_\lambda$
after each pass, depending on the results of the previous pass;
this could be used,
e.g., to emphasize certain locations at fine scales if coarser scale
coefficients
indicate the possible existence of an edge. The results of this paper
guarantee
that each pass will converge.

In this paper we have restricted ourselves to penalty functions that are
weighted $\ell^p$--norms of the $\fg = \left<f, \varphi_{\gamma} \right>$. The
approach can be extended naturally to include penalty functions that
can be written as sums, over $\gamma \in \Gamma$, of more general
functions of $\fg$, so that the functional to be minimized is then written
as
$$
\widetilde{\Phi}_{_{\mbox{\scriptsize{\bf{W}}}}}(f) = \|Kf-g\|^2 +
\sum_{\gamma \in \Gamma}
W_{\gamma} (|\fg|) ~.
$$
The arguments in this paper will still be applicable to this more general case
if the functions $W_{\gamma}: \R_+ \rightarrow \R_+$ are convex, and satisfy
some extra technical conditions, which ensure that the corresponding
generalized component--shrinkage functions $\widetilde{S}_{\gamma}$ are still
non-expansive (used in several places), and that, for some $c > 0$,
$$
\inf_{\|v\| \leq 1} ~ \inf_{\|a\| \leq c} \|v\|^{-2}
\sum_{\gamma \in \Gamma} \left|v_{\gamma} + \widetilde{S}_{\gamma}(a_{\gamma})
-\widetilde{S}_{\gamma}(v_{\gamma}+a_{\gamma}) \right|^2 > 0 ~
$$
(used in Lemma \ref{lm-3-16}).
To ensure that both conditions are satisfied, it is sufficient to choose
functions $W_{\gamma}$ that are convex, with a minimum at $0$ and e.g.
twice differentiable, except possibly at $0$ (where they should nevertheless
still be left and right differentiable), and for which $W_{\gamma}'' >1$
on $V \setminus \{0\}$, where $V$ is a neighborhood of $0$.

We conclude this paper with some comments
concerning the numerical complexity of the algorithm.

At each iteration step,
we must compute the action of the operator $K^*K$ on the current
object estimate, expressed in the $\varphi_{\gamma}$--basis.
In a finite-dimensional setting where the solution is
represented by a vector of length $N$, this necessitates in principle a
matrix multiplication
of complexity $O(N^2)$,
if we neglect the cost of the shrinkage operation in each iteration step.
After sufficient accuracy is attained and the iterations are stopped, the
resulting $(f^n)_{\gamma}$ must be transformed back into the standard
representation domain of the object function, except in the special case
where the $\varphi_{\g}$ are already the basis for the standard representation
(e.g., if the $\varphi_{\g}$ correspond to the pixel representation
for images). This adds one final $O(N^2)$--matrix multiplication. In this
scenario, the total cost equals that
of the classical Landweber algorithm
on the basis of a comparable number of iterations.
Since Landweber's algorithm typically requires a substantial number of
iterations, it follows that this method can become
computationally competitive with the $O(N^3)$ SVD algorithms only
when $N$ is large compared to the number of iterations necessary.

Several methods have been proposed in the literature to accelerate the
convergence
of Landweber's iteration, which could be used for the present algorithm as
well.
For instance, one could use some form of preconditioning (using the
operator $D$
of the Remark \ref{op-D}) or group together $k$ Landweber iteration steps and
apply thresholding only every $k$ steps (see e.g. the book \cite{Eng96}).

Much more substantial gains can be obtained when the operator $\K K$ can be
implemented via fast algorithms. In a first important class of applications,
the matrix \\
$\left( \left<\K K\varphi_{\gamma},\varphi_{\gamma'}\right> \right)_
{\gamma , \gamma' \in \Gamma}$ is sparse; if, for instance, there are only
$O(N)$ non vanishing entries in this matrix, then standard techniques to deal
with the action of sparse matrices will reduce the cost of each iteration step
to $O(N)$ instead of $O(N^2)$. If the $\vpg$--basis
is a wavelet basis, this is the case for a large class
of integro-differential operators of interest (see e.g. \cite{Bey91}).
Even if $\K K$ is sparse in the $\vpg$--basis, but has an even simpler
expression in another basis, and if the transforms back and forth between the
two bases can be carried out via fast algorithms, then it may be useful to
implement the action of $\K K$ via these back--and--forth transformations.
For instance, if the object is of a type that will have a sparse representation
in a wavelet basis, and the operator $\K K$ is a convolution operator, then
we can pick the $\vpg$--basis to be a wavelet basis, and implement the
operator $\K K$ by doing, successively, a fast reconstruction from wavelet
coefficients, followed by a FFT, a multiplication in the Fourier domain, an
inverse FFT, and a wavelet transform, for a total complexity of $O(N \log N)$.
One can obtain similar complexity estimates if the algorithm is modified
to not only take the nonlinear thresholding into account, but also additional
projections $P_{\cal C}$ on a convex set, such as the cone of functions that
are a.e. positive; in this case, after the thresholding operation, one needs
to carry out an additional fast reconstruction from, say, the wavelet domain,
take the positive part, and then perform the fast transform back, without
affecting the $O(N \log N)$ complexity estimate.

The situations described above cover several applications of
great practical relevance, in which we expect this algorithm will prove itself
to be an attractive competitor to other fast techniques for large-scale
inverse problems with sparsity constraints.

\section*{Acknowledgments}

We thank Albert Cohen, Rich Baraniuk, Mario Bertero, Brad Lucier,
St\'ephane Mallat
and especially David Donoho for interesting and stimulating discussions. We
also
would like to thank Rich Baraniuk for drawing our attention to \cite{Fig03}.

Ingrid Daubechies gratefully acknowledges partial support by NSF grants
DMS-0070689
and DMS-0219233, as well as by AFOSR grant F49620-01-1-0099, whereas research
by Christine De Mol is supported by the ``Action de Recherche Concert\'ee'' Nb
02/07-281 and IAP-network in Statistics P5/24.

\section*{Appendices}
\appendix

\section{Wavelets and Besov spaces}
\label{WavBes}

We give a brief review of basic definitions of wavelets and their
connection with Besov spaces.
This will be a sketch only; for details, we direct the reader to e.g.
\cite{Mey92, DeV98, Coh00, Mal98}.

For simplicity we start with dimension 1.  Starting from a (very special)
function
$\psi$ we define\begin{equation*}\psi_{j,k}(x)= 2^{j/2}\ \psi(2^j x-k) ~,
j,k \in \Z~,\end{equation*}
and we assume that the collection $\{\psi_{j,k}; j,~k \in \Z\}$
constitutes an orthonormal basis of $L^2(\mathbb R)$.
For all wavelet bases used in practical applications, there also exists an
associated
{\em scaling function} $\phi$, which is orthogonal to its
translates by integers, and such
that, for all $j \in \Z$,
\begin{equation}
\label{MRA}
\overline{\mbox{Span}\{\phi_{j,k}; k \in \Z\}}~\mbox{\small{$\bigoplus$}}
~\overline
{\mbox{Span}\{\psi_{j,k}; k \in \Z\}}
= \overline{\mbox{Span}\{\phi_{j+1,k}; k \in \Z\}} ~,
\end{equation}
where the $\phi_{j,k}$ are defined analogously to the $\psi_{j,k}$.
Typically, the functions $\phi$ and $\psi$ are very well localized, in the
sense
that $\forall N \in \N$, $\int_{\R} (1+|x|)^N(|\phi(x)|+|\psi(x)|) dx <
\infty$;
one can even
choose $\phi$ and $\psi$ such that they are supported on a finite interval.
This can be
achieved with arbitrary finite smoothness, i.e. for any preassigned $L \in
\N$, one can
find such $\phi$ and $\psi$ that are moreover in $C^L(\R)$. Because of
\eref{MRA},
one can consider (inhomogeneous) wavelet expansions, in which not all
scales $j$
are used,
but a cut-off is introduced at some coarsest scale, often set at $j=0$.
More precisely,
we shall use the following wavelet expansion
of $f
\in L^2$,
\begin{equation}
\label{inhMRA}
 f= \sum_{k=-\infty}^{+\infty}  \left<f,\phi_{0,k}\right> \phi_{0,k} +
\sum_{j=0}^{+\infty}
\sum_{k=-\infty}^{+\infty} \left<f,\psi_{j,k}\right> \psi_{j,k}~.
\end{equation}
Wavelet bases in
higher dimensions can be built by taking appropriate products of
one-dimensional
wavelet and scaling functions. Such $d$-dimensional bases can be viewed as the
result of translating (by elements $k$ of $\Z^d$) and
dilating (by integer powers $j$ of 2) of
not just one, but several
(finite in number) ``mother wavelets'', typically numbered from
1 to $2^d-1$.
It will be convenient to abbreviate the full label (including $j$,  $k$ and the
number of the mother wavelet) to just $\lambda$, with the convention that
$|\lambda|=j$.
We shall again cut off at some coarsest scale, and we shall follow the
convenient
slight abuse
of notation used in \cite{Coh00} that sweeps up the coarsest-$j$
scaling functions
(as in \eref{inhMRA}) into the $\Psi_{\lambda}$ as well. We thus denote the
complete
$d$-dimensional, inhomogeneous wavelet basis by
$\{\Psi_{\lambda}; \lambda \in \Lambda\}$.

It turns out that $\{\Psi_{\lambda}; \lambda \in \Lambda\}$
is not only an orthonormal basis
for $L^2(\R^d)$, but also an unconditional basis for a variety of other
useful
Banach spaces of functions, such
as H\"older spaces, Sobolev spaces and, more generally, Besov spaces.
Again, we review only some basic facts; a full study can be found in e.g.
\cite{Mey92, DeV98, Coh00}.  The Besov spaces $B^s_{p,q}(\R^d)$
consist, basically, of functions that ``have $s$ derivatives in $L^p$'';
the parameter $q$ provides
some additional fine-tuning to the definition of these spaces. The norm
$\|f\|_{_{B^s_{p,q}}}$ in a Besov space $B^s_{p,q}$ is traditionally
defined via the
{\em modulus of continuity} of $f$ in $L^p(\R)$, of which an additional
weighted $L^q$-norm
is then taken, in which the integral is over different scales.
We shall not give its details here; for our purposes it suffices that
this traditional Besov norm is equivalent with a norm that can be computed from
wavelet coefficients. More precisely, let us assume that the original
1-dimensional $\phi$ and
$\psi$ are in $C^L(\R)$, with $L>s$, that
$\sigma=s+d(\frac{1}{2}-\frac{1}{p}) \geq 0$,
and define the norm $\VVert \cdot \VVert_{_{s;p,q}}$ by
\begin{equation}
\label{triple}
\VVert f \VVert  _{_{s;p,q} }= \left( \sum_{j=0}^{\infty} \left(2^{j \sigma
p} \sum_{
\lambda \in \Lambda ,
|\lambda |=j }|\left<f,\Psi_{\lambda}\right>|^p\right)^{q/p}\right)^{1/q} ~~.
\end{equation}
Then this norm is equivalent to the traditional Besov norm,
$\VVert f \Vvert _{s;p,q}\sim\| f \|_{_{B^s_{p,q}}}$, that is, there exist
strictly positive constants $A$ and $B$ such that
\begin{equation}
\label{Besnor}
A \VVert f \VVert_{_{s;p,q}} \leq \| f \| _{_{B^s_{p,q}}}
\leq B \VVert f \Vvert_{_{s;p,q}} ~.
\end{equation}
The condition that $\sigma \geq 0$ is imposed to ensure
that $B^s_{p,q}(\R^d)$ is a subspace
of $L^2(\R^d)$; we shall restrict ourselves to this case in this paper.
From \eref{triple} one can gauge the fine-tuning role played by the
parameter $q$
in the definition of the Besov spaces. A particularly convenient choice, to
which we
shall stick in the remainder of this paper, is $q=p$,
for which the expression simplifies
to
\begin{equation*}
\VVert f \VVert_{_{s,p}} = \left( \sum_{\lambda \in \Lambda} 2^{\sigma
p|\lambda|} ~
| \left< f, \Psi_{\lambda} \right> |^p \right)^{1/p} ~~;
\end{equation*}
to alleviate notation, we shall drop the extra index $q$ wherever it
normally occurs,
on the understanding that $q=p$ when we do so.

When $0<p,~q<1$, the Besov spaces can still be defined as complete metric spaces,
although they are no longer Banach spaces (because (\ref{triple}) no longer is a norm),
This allows for more local variability in local smoothness than
is typical for functions
in  the usual H\"older
or Sobolev spaces. For instance, a real function $f$ on $\R$ that is piecewise
continuous, but for which each piece is locally in $C^s$, can be an element of
$B^s_p(\R)$, despite the possibility of discontinuities at the
transition from one piece to
the next, provided $p>0$ is sufficiently small, and
some technical conditions are met on the number and
size of the discontinuities, and on the decay at $\infty$ of $f$. 

Wavelet bases are thus closely linked to a rich class
of smoothness spaces; they also provide a good tool for high accuracy
nonlinear approximation of a wide variety of functions.  
For instance, if the bounded function $f$
on $[0,1]$ has only finitely many discontinuities, and is $C^s$ elsewhere,
then one can find a way of renumbering (dependent on $f$ itself)
the wavelets in the standard wavelet
expansion of $f$, so that the distance in, say, $L^2([0,1])$ between $f$
and the first $N$ terms of this reordered wavelet expansion, decreases
proportionally to $N^{-s}$.
If $s$ is large, it follows that  a very accurate approximation to $f$ can be
obtained with relatively few wavelets; this
is possible because
the smooth patches of the piecewise continuous
$f$  will be well approximated by coarse scale wavelets, which are 
few in number; to capture the behavior of $f$ near 
the discontinuities  much more localized
finer scale wavelets are required, but only 
those wavelets located exactly near
the discontinuities will be needed, which amounts
again to a small number. 

In higher dimensions, $d > 1$,  the suitability of wavelets
is influenced by the
dimension of the manifolds on which singularities occur. If the
singularities in the
functions of interest are solely point singularities, then expansions
using $N$ wavelets can still approximate such functions with distances that
decrease like $N^{-s}$,
depending on their behavior away
from the
singularities. If, however, we are interested in $f$ that may have, e.g.
discontinuities
along manifolds of dimension higher than 0, then such wavelet
approximations  are not optimal.
For instance, if $f:\R^2 \rightarrow \R$ is piecewise $C^L$, with
possible jumps across the
boundaries of the smoothness domains, which are themselves smooth
(say, $C^L$ again) curves,
then $N$-term wavelet approximations  to $f$ cannot achieve an error
rate decay faster than $N^{-1/2}$,
regardless of the
value of $L>1$. 

It follows that whenever we are faced with an inverse problem
that needs regularization,
in which the objects to be restored are expected to be mostly smooth,
with very localized
lower dimensional areas of singularities, 
we can expect that their expansions into wavelets
will be sparse. This sparsity can be expressed by requiring that
the wavelet coefficients (possibly with some scale-dependent weight)
have a finite (or small) $\ell^p$-norm,
with $1\leq p \leq 2$, or equivalently that  the Besov-equivalent norm $\VVert f
\VVert_{_{s,p}}$ is finite (or small), where $\VVert f \VVert_{_{s,p}}$
is exactly of the form
$\Phi_{\mathbf w,p}$ defined in \eref{funct-gen}.
 
\section{ A fixed-point theorem}
\label{Opial}

We provide here the proof of the theorem needed to establish the weak
convergence of the iterative algorithm. The theorem is given in \cite{Opi67};
we give a simplified proof here (see the remark at the end),
which nevertheless still follows the
main lines of Opial's paper.

\begin{theorem}
\label{FPThm}
Let ${\cal C}$ be a closed convex subset of the Hilbert space $\cH$ and let
the mapping $\A : {\cal C} \to {\cal C}$ satisfy the following conditions:
\begin{enumerate}
\item[{\rm (i)}] $\A $ is non-expansive: $\| \A  v -  \A
 v' \| \leq \| v - v'\|,\ \forall v,v' \in {\cal C}$~,
\item[\rm{ (ii)}] $\A $ is asymptotically regular: $\| \A ^{n+1}
v -\A ^n v\| \xrightarrow[n \to \infty]{~} 0,\ \forall v \in
{\cal C}$~,
\item[\rm{ (iii)}] the set ${\cal F}$ of the fixed points of $\A $ in
${\cal C}$ is not empty~.
\end{enumerate}
Then, $\forall v \in \cal C$, the sequence $(\A ^n v)_{n \in \mathbb{N}}$
converges weakly to a fixed point in ${\cal F}$.
\end{theorem}
The proof of the main theorem will follow from a series of lemmas.
As before, we use the notation {\em w}$\,$--$\lim$ to indicate a {\em weak}
limit.
\begin{lemma}
\label{FP1}
If $u,v \in\cH$, and if $(v_n)_{n \in \N}$ is a sequence in $\cH$ such that
w--$\lim_{n \to \infty}v_n = v$, and $u \neq
v$, then
$\lim\inf_{n \to \infty} \| v_n - u\| > \lim\inf_{n \to \infty}\| v_n - v\|~$.
\end{lemma}
{\em Proof:}
We have $\lim\inf_{n \to \infty}\| v_n - u\|^2 $
$= \lim\inf_{n \to \infty}\| v_n - v\|^2  +
\| v - u\|^2 +  2 \lim_{n \to \infty} Re (v_n-v,v-u)$
$= \lim\inf_{n \to \infty}\| v_n - v\|^2  +\| v - u\|^2~$,
whence the result.
\hfill\QED

\bigskip

\begin{lemma}
\label{FP2}
Suppose that $\A:\cal C \rightarrow \cal C$ satisfies condition
{\rm(i)} in Theorem \ref{FPThm}.\\
If  w--$\lim_{n \to \infty}u_n= u$, and
$\lim_{n \to \infty} \|u_n -\A u_n -h\| =0~$, then
$h = u -\A u$\ .
\end{lemma}
{\em Proof:}
Because of the non-expansivity of $\A $ (assumption (i)),
we have
$\| u_n - (h + \A u)\| $
$\leq \| u_n - h - \A  u_n\|$ $+\|  \A  u_n - \A  u\| $
$\leq \| u_n - h - \A  u_n\|$ $ +\|   u_n -  u \|~$.
Hence,
\begin{eqnarray*}
{\mathop{\rm lim\ inf}_{n \to \infty}}\ \| u_n - (h + \A  u)\|
&\leq&
{\mathop{\rm lim}_{n \to \infty}}\ \| h - (u_n - \A  u_n)\| +
{\mathop{\rm lim\ inf}_{n \to \infty}}\ \| u_n -  u\| \\
&=&
{\mathop{\rm lim\ inf}_{n \to \infty}}\ \| u_n -  u\|
\end{eqnarray*}
It then follows from Lemma \ref{FP1} that $u = h + \A  u$ or $h = u -
\A  u$.
\hfill\QED

\bigskip

\begin{lemma}
\label{FP3}
Suppose that $\A :\,\cal C \rightarrow \cal C$ satisfies conditions
{\rm(i)} and {\rm(ii)} in Theorem \ref{FPThm}.
If a subsequence of $(\A ^n v)_{n\in \mathbb{N}}$, with $v \in {\cal C}$,
converges weakly in ${\cal C}$, then its limit is in $\cal F$.
\end{lemma}
{\em Proof:}
Suppose {\em w}$\,$--$\lim_{k \to \infty}\A ^{n_k}v= u~$.
Since, by the assumption (ii) of asymptotic
regularity,
$\lim_{n \to \infty}\|\A ^{n} v - \A  \A ^{n} v \|=0~$,
we have
$\lim_{k \to \infty}\|\A ^{n_k} v - \A  \A ^{n_k} v \|=0~$.
By Lemma \ref{FP2}, it
follows that $u - \A u = 0$ , i.e. that $u$ is in $\cal F$.
\hfill\QED

\bigskip

\begin{lemma}
\label{FP4}
Suppose that $\A :\,\cal C \rightarrow \cal C$ satisfies conditions
{\rm(i)} and {\rm(iii)} in Theorem \ref{FPThm}. Then,
for all $h \in {\cal F}$, and all $v \in {\cal C}$, the sequence
$(\|\A ^n v -
h\|)_{n\in \N}\ $ is non-increasing and thus has a limit.
\end{lemma}
{\em Proof:}
Since $\A $ is non-expansive, we have indeed
$\| \A ^{n+1} v - h \|$
$= \| \A \A ^n v -\A   h\| $
$\leq \| \A ^{n} v - h \|~.$
$~~~~~~~~~~~~~~~~~~~~~~~~$\hfill\QED

\bigskip

We can now proceed to the

\bigskip

\noindent
{\bf Proof of Theorem \ref{FPThm}}

\noindent
Let $v$ be any element in $\cal{C}$. Take an arbitrary $h \in {\cal F}$.
By Lemma \ref{FP4}, we then have \\
$\lim\sup_{n \to \infty}\| \A ^n v\|$
$ \leq
\lim\sup_{n \to \infty} \|\A ^n v - h\|$
$+\| h\|$
$ = \| h\| $  $+ \lim_{n \to \infty}\ \|\A ^n v - h\|$
$< \infty~$.

\noindent
Since the $\|\A ^n v\|$ are thus uniformly bounded,
it follows from the Banach-Alaoglu theorem that they must have at least
one weak accumulation point.

\noindent
The following argument shows that this
accumulation point is unique.
Suppose we have two different accumulation points :
{\em w}$\,$--$\lim_{k \to \infty}\A ^{n_k} v =u$, and
{\em w}$\,$--$\lim_{\ell \to \infty}\A ^{{\tilde n}_\ell} v =\tilde{u}~$,
with $u \neq {\tilde u}$.

\noindent
By Lemma \ref{FP3}, $u$ and $\tilde u$ must both lie in $\cal F$,
and by Lemma \ref{FP4},
the limits $\lim_{n \to \infty} \|\A ^n v -
u\|$ and $\lim_{n \to \infty} \|\A ^n v -
{\tilde u} \|$ both exist.

\noindent
Since $\tilde{u} \neq u$, we obtain from Lemma \ref{FP1} that
$\lim\inf_{k \to \infty} \|\A ^{n_k} v - {\tilde u}\| $
$ > {\lim\inf}_{k \to \infty} \|\A ^{n_k} v - u\|\ .$
On the other hand,
because $(\|\A ^{n_k} v-{\tilde u}\|)_{k\in \mathbb{N}}$ and
$(\|\A ^{n_k} v-u\|)_{k\in \mathbb{N}}$
are each a
subsequence of a convergent sequence,
$\lim\inf_{k \to \infty} \|\A ^{n_k} v - {\tilde u}\|$ =
$\lim_{n \to \infty} \|\A ^n v - {\tilde u}\|$ and
$\lim\inf_{k \to \infty} \|\A ^{n_k} v - u\|$ =
$ \lim_{k \to \infty} \|\A ^{n_k} v - u\|~$.
It follows that
$\lim_{n \to \infty} \|\A ^{n} v - {\tilde
u}\|$ $> \lim_{n \to \infty} \|\A ^n v -u\| ~$.
In a completely analogous way (working with the subsequence
$\A ^{{\tilde n}_l} v$ instead of $\A ^{n_k} v$) one derives
the opposite strict inequality. Since both cannot be valid
simultaneously, the assumption of the existence of two different
weak accumulation points for $(\A ^n v)_{n\in \N}$ is false.

\noindent
It thus follows that $\A ^n v$ converges weakly to this unique
weak accumulation point.
\hfill\QED

\bigskip

\begin{remark}
{\rm It is essential to require that the set $\cal F$ is not empty since
there are
asymptotically regular, non-expansive maps that possess no fixed point.
However, the only place where we used this assumption was in showing that the
$\|\A ^n v \|$ were bounded. If one can prove this boundedness by
some other means (e.g. by a variational principle as we did in the iterative
algorithm), then we automatically have a weakly convergent subsequence
$(\|\A ^{n_k} v\|)_{k\in \mathbb{N}}$, and thus, by Lemma
\ref{FP3}, an element of $\cal F$.}
\end{remark}

\begin{remark}
{\rm The simplification of the original argument of \cite{Opi67}
(obtained through
deriving the contradiction in the proof of Theorem \ref{FPThm}) avoids having
to appeal to the convexity of $\cal F$ (which is true but not immediately
obvious) and having to introduce the auxiliary sets $\cal F_\delta$ used in
\cite{Opi67}.}
\end{remark}


\begin{thebibliography}{A}

\bibitem [1]{Abr98} F. Abramovich and B. W. Silverman, \textit{Wavelet
Decomposition Approaches to Statistical Inverse Problems.} Biometrika
\textbf{85} (1998), 115--129.

\bibitem [2]{Ber98} M. Bertero and P. Boccacci, \textit{Introduction to
Inverse Problems in Imaging}, Institute of Physics, Bristol, 1998.

\bibitem [3]{Ber96} M. Bertero and C. De Mol, \textit{Super-resolution by
data inversion}, in: Progress in Optics (Vol. XXXVI), E. Wolf, ed.,
Elsevier, Amsterdam, 1996, pp. 129--178.

\bibitem [4]{Bey91} G. Beylkin, R. Coifman and V. Rokhlin, \textit{Fast
Wavelet Transforms and Numerical Algorithms I.} Comm. Pure Appl. Math.
\textbf{44} (1991), 141--183.

\bibitem [5]{Can00} E. J. Cand\`es and D. L. Donoho, \textit{Recovering Edges
in Ill-Posed Inverse Problems: Optimality of Curvelet Frames.} Ann. Statist.
\textbf{30} (2000), 784--842.

\bibitem [6]{Cha98} A. Chambolle, R. A. DeVore, N.-Y. Lee and B. J. Lucier,
\textit{Nonlinear Wavelet Image Processing: Variational Problems, Compression,
and Noise Removal through Wavelet Shrinkage.} IEEE Trans. Image Processing
\textbf{7} (1998), 319--335.

\bibitem [7]{CDS01} S. Chen, D. Donoho and M. Saunders, \textit{Atomic
Decomposition by Basis Pursuit} SIAM Review \textbf{43} (2001), 129--159.

\bibitem [8]{Coh00} A. Cohen, \textit{Wavelet methods in numerical
analysis.}, Handbook of Numerical Analysis, vol. VII, P. G. Ciarlet and J. L.
Lions eds., Elsevier, Amsterdam, 2000.

\bibitem [9]{Coh02} A. Cohen, M. Hoffmann and M. Reiss, \textit{Adaptive
wavelet Galerkin methods for linear inverse problems.} preprint, 2002.


\bibitem [10]{DeM02} C. De Mol and M. Defrise, \textit{A note on
wavelet-based inversion methods}, in: \textit{Inverse Problems, Image Analysis
and Medical Imaging},  M. Z. Nashed and O. Scherzer eds,
Series ``Contemporary Mathematics''
vol. 313, pp. 85--96, American Mathematical Society, 2002.

\bibitem [11]{DeP95} A. R. De Pierro, \textit{A modified expectation
maximization algorithm for penalized likelihood estimation in emission
tomography.} IEEE Trans. Med. Imag. \textbf{14} (1995), 132--137.

\bibitem [12]{DeV98} R. DeVore, \textit{Nonlinear Approximation.} Acta
Numerica (1998), 1--99.

\bibitem [13]{Dic96} V. Dicken and P. Maass, \textit{Wavelet-Galerkin
methods for ill-posed problems.} J. Inv. Ill-Posed Problems  \textbf{4} (1996)
203--222.

\bibitem [14]{Don92} D. Donoho, \textit{Superresolution via sparsity
constraints.}
SIAM J. Math. Anal. \textbf{23} (1992), 1309--1331.

\bibitem [15]{Don95} D. Donoho, \textit{Nonlinear solution of Linear
Inverse Problems by Wavelet-Vaguelette Decomposition.} Appl. Comp. Harmonic
Anal. \textbf{2} (1995), 101--126.

\bibitem [16]{Don00} D. Donoho, \textit{Orthonormal ridgelets and linear
singularities.} SIAM J. Math. Anal.
\textbf{31} (2000), 1062--1099.

\bibitem [17]{Don94} D. Donoho and I. Johnstone, \textit{Ideal
spatial adaptation via wavelet shrinkage.} Biometrika
\textbf{81} (1994), 425--455.

\bibitem [18]{DuSh52} R. J. Duffin and A. C. Schaeffer, \textit{A class of
nonharmonic Fourier series.} Trans. Am. Math. Soc. \textbf{72} (1952),
341--366.

\bibitem [19]{Eic92} B. Eicke, \textit{Iteration methods for convexly
constrained ill-posed problems in Hilbert space.} Numer. Funct. Anal. Optim.
\textbf{13} (1990), 413--429.

\bibitem [20]{Eng96} H. W. Engl, M. Hanke and A. Neubauer,
\textit{Regularization of Inverse Problems}, Kluwer, Dordrecht, 1996.

\bibitem [21]{Fig03} M. Figueiredo and R. Nowak, \textit{An EM Algorithm
for Wavelet-Based Image Restoration}, IEEE Transactions on Image Processing.
To appear in July 2003.

\bibitem [22]{KMR03} J. Kalifa, S. Mallat and B. Roug\'e,
\textit{Deconvolution by thresholding in mirror wavelet bases.} IEEE Trans. on
Image Processing \textbf{12} (2003), 446--457.

\bibitem [23]{Lad51} L. Landweber, \textit{An iterative formula for Fredholm
integral equations of the first kind.} Am. J. Math.
\textbf{73} (1951), 615--624.

\bibitem [24]{Lan00}  K. Lange, D. R. Hunter and I. Yang,
\textit{Optimization Transfer algorithms using surrogate objective functions.}
J. Comp. Graph. Stat. \textbf{9} (2000), 1--59.

\bibitem [25]{Lee01}  N.-Y. Lee and B. J. Lucier,
\textit{Wavelet Methods for Inverting the Radon Transform with Noisy Data.}
IEEE Trans. Image Processing \textbf{10} (2001), 79--94.

\bibitem [26]{Li02} M. Li, H. Yang and H. Kudo, \textit{An accurate iterative
reconstruction algorithm for sparse objects: application to 3-D blood-vessel
reconstruction from a limited number of projections.} Phys. Med. Biol
\textbf{47}
(2002), 2599--2609.

\bibitem [27]{Lou97} A. K. Louis, P. Maass and A. Rieder, \textit{Wavelets:
Theory and Applications}, Wiley, Chichester, 1997.

\bibitem [28]{Mal98} S. Mallat, \textit{A Wavelet Tour of Signal
Processing}, 2nd edition, Academic Press, San Diego, 1999.

\bibitem [29]{Mey92} Y. Meyer, \textit{Wavelets and Operators}, Cambridge
University Press, 1992.

\bibitem [30]{Nov01} R. Nowak and M. Figueiredo, \textit{Fast wavelet-based
image deconvolution using the EM algorithm.} Conference Record of the
Thirty-Fifth
Asilomar Conference on Signals, Systems and Computers, Vol. 1 ,
pp. 371--375, 2001.

\bibitem [31]{Opi67} Z. Opial, \textit{Weak convergence of the sequence of
successive approximations for nonexpansive mappings.} Bull. Amer. Math. Soc.
\textbf{73} (1967), 591--597.

\bibitem [32]{Tib96} R. Tibshirani, \textit{Regression shrinkage and
selection via the lasso.} J. Royal Statist. Soc. B \textbf{58} (1996),
267--288.


\end{thebibliography}
\end{document}